\documentclass[11pt,reqno]{amsart}

\usepackage[latin1]{inputenc} 
\usepackage{amsthm,amsmath,amssymb,bm,mathtools,setspace} 
\usepackage{graphicx}
\usepackage{color}
\usepackage{verbatim}
\usepackage{float}
\usepackage{nccmath}
\usepackage[noadjust]{cite}
\usepackage[foot]{amsaddr}
\usepackage{enumitem}
\usepackage{thmtools}
\usepackage{thm-restate}
\usepackage{eqparbox}

\setlist[enumerate,1]{label=\upshape(\roman*),ref=(\roman*)}

\newcommand{\COMMENT}[1]{}
\renewcommand{\COMMENT}[1]{\footnote{\textcolor{blue!70!black}{#1}}} 

\usepackage{hyperref}       
\usepackage[capitalise]{cleveref} 
\hypersetup{colorlinks=true,
   citecolor=blue,
   filecolor=blue,
   linkcolor=blue,
   urlcolor=blue
}

\creflabelformat{equation}{(#2#1#3)}
\crefdefaultlabelformat{#2#1#3}

\Crefname{enumi}{}{}
\Crefname{thm}{Theorem}{Theorems}
\Crefname{lm}{Lemma}{Lemmas}
\Crefname{cor}{Corollary}{Corollaries}
\Crefname{prop}{Proposition}{Propositions}
\Crefname{claim}{Claim}{Claims}
\Crefname{equation}{}{}
\Crefname{conjecture}{Conjecture}{Conjectures}
\Crefname{figure}{Figure}{Figures}
\Crefname{fact}{Fact}{Facts}

\usepackage{array}  
\newcolumntype{C}{>{\centering \arraybackslash}m{3cm}}
\newcolumntype{D}{>{\centering \arraybackslash}m{1.25cm}}

\usepackage{enumitem} 
\setlist[enumerate]{topsep=0pt,itemsep=-0.5ex,partopsep=1ex,parsep=1ex}
\setlist[itemize]{topsep=0pt,itemsep=-0.5ex,partopsep=1ex,parsep=1ex}

\usepackage{bbm}

\theoremstyle{plain}
\newtheorem{theorem}{Theorem}[section]

\newtheorem{lemma}[theorem]{Lemma}
\newtheorem{cor}[theorem]{Corollary}

\newtheorem{conj}[theorem]{Conjecture}

\theoremstyle{definition}

\renewcommand{\epsilon}{\varepsilon}

\newcommand{\stri}{\scriptscriptstyle\triangle}
\newcommand{\splus}{\scriptscriptstyle+}

\def \Gs {G_{\mathrm{s}}}
\def \Gc {G_{\mathrm{c}}}

\newcommand{\D}{\mathcal{D}}
\newcommand{\G}{\mathcal{G}}
\newcommand{\cD}{\mathcal{D}}

\newcommand{\T}{\mathcal{T}}

\renewcommand{\P}{\mathcal{P}}

\renewcommand{\leq}{\leqslant}
\renewcommand{\geq}{\geqslant}

\usepackage[dvipsnames]{xcolor}

\usepackage{stackengine}
\stackMath
\newcommand{\stacklEq}[1]{%
  \setbox0=\hbox{${}\mathrel{\stackon[-1pt]{\leq}{\scriptstyle\text{#1\strut}}}{}$}
  \xdef\tmpwd{\dimexpr\the\wd0\relax}
  \kern.5\tmpwd\mathclap{\box0}&\kern.5\tmpwd
}

\newcommand{\stackEq}[1]{
  \setbox0=\hbox{${}\mathrel{\stackon[-1pt]{=}{\scriptstyle\text{#1\strut}}}{}$}
  \xdef\tmpwd{\dimexpr\the\wd0\relax}
  \kern.5\tmpwd\mathclap{\box0}&\kern.5\tmpwd
}

\allowdisplaybreaks

\author[B. Granet]{Bertille Granet}
\email{first.lastname [at] warwick.ac.uk}
\address[Granet]{Warwick Mathematics Institute, University of Warwick, Coventry CV4 7AL, United Kingdom.}
\author[D. Horsley]{Daniel Horsley}
\email{first.lastname [at] monash.edu}
\address[Horsley]{School of Mathematics, Monash University, Melbourne VIC 3800, Australia.}
\thanks{This work was partially supported by the Australian Research Council grant DP240101048. Bertille Granet was also partially supported by a Leverhulme Early Career Fellowship.}

\title{Determining decomposition thresholds for long odd cycles}

\begin{document}

\begin{abstract}
An \emph{$\ell$-cycle decomposition} of a graph $G$ is a set of $\ell$-cycles in $G$ whose edge sets partition the edge set of $G$. The \emph{$\ell$-cycle decomposition threshold} $\delta_{C_\ell}$ is then the least real number such that any $n$-vertex graph $G$ with minimum degree at least $(\delta_{C_\ell}+o(1))n$ has an $\ell$-cycle decomposition if and only if $\ell$ divides $|E(G)|$ and each vertex of $G$ has even degree.

Nash-Williams' famous conjecture on triangle decompositions states, asymptotically, that $\delta_{C_3}=\frac{3}{4}$. A very recent breakthrough result of Delcourt and Postle completely resolved this conjecture, however, Glock, K\"uhn, and Osthus have posed the problem of determining $\delta_{C_\ell}$ for larger odd values of $\ell$ (the behaviour of $\delta_{C_\ell}$ for even $\ell$ is different and well understood). A natural generalisation of Nash-Williams' conjecture implies that $\delta_{C_\ell}=\frac{\ell}{2\ell-2}$ for all odd $\ell \geq 3$. Here we prove that this conjecture holds for all $\ell \geq 73$.
\end{abstract}

\maketitle
\setstretch{1.1}

\section{Introduction}

The study of graph decompositions has a very long history. Already back in 1847, Kirkman \cite{Kir} showed that the edges of the complete graph $K_n$ can be partitioned into triangles if and only if $n$ is congruent to $1$ or $3$ modulo $6$. In 1883,
Walecki \cite{Luc} (see also \cite{alspach2008wonderful}) showed that $K_n$ can be decomposed into Hamilton cycles if and only if $n$ is odd. 
More generally, given a graph $F$, one can ask for which values of $n$ there is a decomposition of $K_n$ into copies of $F$ (an \emph{$F$-decomposition}). Wilson \cite{Wil} proved in 1975 that, for large $n$, $K_n$ admits an $F$-decomposition if an only if it is \emph{$F$-divisible}, where a graph $G$ is said to be $F$-divisible if $e(F)$ divides $e(G)$ and $\gcd(F)$ divides $\gcd(G)$, with $\gcd(H)$ denoting the largest integer which divides all the degrees in $H$.

More recently, research has expanded to finding necessary and sufficient conditions for decomposing non-complete graphs, and in particular, as is often the case in modern extremal graph theory, finding optimal minimum-degree conditions which guarantee the existence of a desired decomposition. For instance, Csaba, K\"uhn, Lo, Osthus, and Treglown \cite{CsaKuhLoOstTre} showed that any large $n$-vertex even-regular graph of degree at least $\frac{n}{2}$ has a decomposition into Hamilton cycles (and this is tight as there exist non-Hamiltonian $(\frac{n}{2}-1)$-regular graphs on $n$ vertices). The key conjecture in the area is due to Nash-Williams and states that any sufficiently large $n$-vertex $K_3$-divisible graph of minimum degree at least $\frac{3n}{4}$ has a triangle decomposition. This conjecture attracted a lot of interest over the last few decades (see below for details). In the final stages of the preparation of this manuscript, a fantastic breakthrough was announced by Delcourt and Postle \cite{DelPos2026}, who completely resolved this conjecture.

\begin{theorem}[Delcourt and Postle {\cite{DelPos2026}}]\label{conj:NW}
    Let $n$ be sufficiently large and $G$ be an $n$-vertex graph of minimum degree $\delta(G)\geq \frac{3n}{4}$. If $G$ is $K_3$-divisible, then it admits a $K_3$-decomposition. 
\end{theorem}

This is tight for the following reason. Let $\G_3$ be the set of graphs $G$ on $n$ vertices such that $n \equiv 0 \pmod 4$ and $G$ is obtained from the complete balanced bipartite graph with parts $A$ and $B$ of size $\frac{n}{2}$ by adding an $\frac{n}{4}$-regular graph in each part. To cover the $\frac{n^2}{4}$ edges between $A$ and $B$, one requires $\frac{n^2}{8}$ triangles, each using exactly one of the $\frac{n^2}{8}$ edges within $A$ or $B$. Thus, no graph obtained from a $G \in \G_3$ by removing edges within $A$ or $B$ can have a triangle decomposition. In particular, when for instance $n=12k$ for some odd natural number $k$ and a perfect matching is removed within each of $A$ and $B$, it is not too difficult to verify that the resulting graph is $K_3$-divisible, thus confirming that minimum degree $\frac{3n}{4}$ is indeed required in \cref{conj:NW}.

Obtaining precise minimum-degree bounds which guarantee the existence of an $F$-decomposition is often difficult and has, so far, only been achieved for a few graphs $F$ (see e.g.\ \cite{yuster2000packing} for trees, \cite{Tay} for cycles of even length not equal to 6, \cite{CsaKuhLoOstTre} for Hamilton cycles, and now \cite{DelPos2026} for triangles).
This gives rise to the notion of \emph{decomposition threshold}, which instead captures the \emph{asymptotically} optimal minimum-degree condition required.
More precisely, define the \emph{decomposition threshold} $\delta_F$ of a graph $F$ to be 
the least $\delta$ such that, for all $\varepsilon >0$, there exists $n_0$ such that any $F$-divisible graph on $n\geq n_0$ vertices with minimum degree at least $(\delta+\varepsilon)n$ has an $F$-decomposition.

Finding decomposition thresholds is still very hard, and, to our knowledge, they have so far only been explicitly determined for bipartite graphs $F$ \cite{GloKuhLoMonOst} and very recently for triangles \cite{DelPos2026}.
To determine or bound $\delta_F$, other weaker parameters have been of crucial importance. We now proceed to define them.
Let $G$ be a graph on $n$ vertices. Given a graph $F$ and $\eta>0$, an \emph{$\eta$-approximate $F$-decomposition} of $G$ is a collection of edge disjoint copies of $F$ which cover all but at most $\eta n^2$ edges of $G$; and a \emph{fractional $F$-decomposition} of $G$ is a function $\varphi$
which assigns a real number in $[0,1]$ to each copy of $F$ in $G$ and such that for each edge $e\in E(G)$, the sum of the values of $\varphi$ over the copies of $F$ in $G$ that contain $e$ is equal to $1$. 
One can then define \emph{approximate} and \emph{fractional} analogues of the decomposition threshold.
Namely, for $\eta>0$, define the \emph{$\eta$-approximate decomposition threshold} for $F$ to be the 
least $\delta$ such that, for all $\varepsilon>0$, there exists $n_0$ such that any $F$-divisible graph on $n\geq n_0$ vertices with minimum degree at least $(\delta+\varepsilon)n$ has an $\eta$-approximate $F$-decomposition. The \emph{approximate decomposition threshold} is then $\delta_F^{0+}\coloneqq \sup_{\eta>0}\delta_F^\eta$.
Similarly, let $\delta_F^*$ be the \emph{fractional decomposition threshold} for $F$, defined as 
the least $\delta$ such that, for all $\varepsilon>0$, there exists $n_0$ such that any $F$-divisible graph on $n\geq n_0$ vertices with minimum degree at least $(\delta+\varepsilon)n$ has a fractional $F$-decomposition.
Observe that $F$-divisibility is of course not necessary for the existence of an approximate or fractional $F$-decomposition, but this technical condition is traditionally added to these threshold definitions to ensure that $\delta_F^{0+},\delta_F^*\leq \delta_F$.

The decomposition threshold $\delta_F$ has strong ties with its approximate and fractional counterparts $\delta_F^{0+},\delta_F^*$. Indeed, Haxell and R\"odl \cite{HaxRod} first showed using the regularity lemma that any fractional decomposition can be turned into an approximate decomposition and hence $\delta_F^{0+}\leq \delta_F^*$. Then, using the iterative absorption method, Barber, K\"uhn, Lo, and Osthus \cite{BarKuhLoOst}, and Glock, K\"uhn, Lo, Montgomery, and Osthus \cite{GloKuhLoMonOst} proceeded to show that, for certain graphs $F$ such as odd cycles 
and cliques (respectively), the existence of approximate $F$-decompositions above a certain minimum degree threshold implies the existence of full decompositions above that same threshold and hence $\delta_F= \delta_F^{0+}$.

\begin{theorem}[Haxell and R\"odl \cite{HaxRod}, Barber, K\"uhn, Lo, and Osthus \cite{BarKuhLoOst}, Glock, K\"uhn, Lo, Montgomery, and Osthus \cite{GloKuhLoMonOst}]\label{T:reduction}
    Let $F$ be an odd cycle 
    or a clique. Then, $\delta_F=\delta_F^{0+}=\delta_F^*$.
\end{theorem}

This result implies that, to determine the decomposition threshold of odd cycles
and cliques, it suffices to consider its fractional counterpart. This has been of crucial importance towards the proof of \cref{conj:NW}.
Following multiple preceding results \cite{Gar,Dro,DukHor}, the best fractional bound for triangles, due to Delcourt and Postle \cite{DelPos2021} was previously $\delta_{K_3}^*\leq \frac{7+\sqrt{21}}{14}\leq 0.82733$. 
In their very recent breakthrough \cite{DelPos2026}, the same authors show that indeed $\delta_{K_3}^*=\frac{3n}{4}$, which together with \cref{T:reduction} already gives an asymptotic version of \cref{conj:NW}. Through a careful analysis of extremal examples and the use of their recent refined absorption method, they proceed to prove a stability theorem, which leads to the full resolution of Nash-Williams' conjecture (\cref{conj:NW}).

For cliques, a folklore generalisation of Nash-Williams' conjecture stated that $\delta_{K_r}=1-\frac{1}{r+1}$ for all $r\geq 4$. Surprisingly, this was recently disproved by Delcourt, Henderson, Lesgourgues, and Postle \cite{DelHenLesPos}, who showed that for all $r\geq 4$, there exists $c>1$ such that $\delta_{K_r}>1-\frac{1}{c(r+1)}$. The current best upper bounds are $\delta_{K_4}\leq 1-\frac{2}{33}$ due to Zhang and Ge \cite{ZhaGe} and $\delta_{K_r}\leq 1-\frac{1}{100 r}$ for all $r\geq 5$ due to Montgomery \cite{montgomery2019fractional}.

The decomposition threshold of bipartite graphs $F$ has been fully determined. Indeed, Yuster \cite{Yus} first showed that $\delta_F=\frac12$ when $F$ is a connected bipartite graph with a vertex of degree $1$, while Barber, K\"uhn, Lo, and Osthus \cite{BarKuhLoOst} showed that $\delta_{C_4}=\frac{2}{3}$ and $\delta_{C_\ell}=\frac12$ for any even $\ell\geq 6$. The value of $\delta_F$ for other bipartite graphs $F$ was then fully resolved by Glock, K\"uhn, Lo, Montgomery, and Osthus \cite{GloKuhLoMonOst}.

Note that the works \cite{BarKuhLoOst,GloKuhLoMonOst} also make use of the iterative absorption method to determine $\delta_F$ for bipartite $F$.
For non-bipartite graphs $F$, this method does not immediately give an exact value for $\delta_F$ as the minimum degree required to find suitable absorbers cannot always be shown to correspond to the approximate decomposition threshold, leading to more complicated functions for $\delta_F$ in terms of $\delta_F^{0+}$ or $\delta_F^*$. For instance, for graphs $F$ of chromatic number $\chi(F)\geq 5$, Glock, K\"uhn, Lo, Montgomery, and Osthus \cite{GloKuhLoMonOst} showed that $\delta_F\in \{\delta_F^*, 1-\frac{1}{\chi(F)}, 1-\frac{1}{\chi(F)+1}\}$ and conjectured that none of the terms in this discretisation result can be omitted.
For graphs of chromatic number $3$ or $4$, they do not obtain an analogous discretisation result and expect a different behaviour. Together with previous bounds of \cite{GloKuhLoMonOst}, the main result of \cite{DelPos2026} implies that $\delta_F^*\leq \frac{3}{4}$ for any graph $F$ of chromatic number $3$.
To investigate the cases $\chi(F)\in \{3,4\}$ further, the problem of determining the decomposition thresholds of odd cycles was posed by Glock, K\"uhn, and Osthus in their survey \cite{GloKuhOst}, which we direct the reader to for an in-depth overview of further results on decomposition thresholds.

In this paper, we solve all but a small number of cases of this problem and determine $\delta_{C_\ell}$ for any odd $\ell\geq 73$.
A natural lower bound can be obtained by generalising mutatis mutandis the construction of the extremal  
graphs in $\G_3$ to longer odd cycles. Indeed, let $\G_\ell$ be the set of graphs $G$ on $n$ vertices such that $n \equiv 0 \pmod{2\ell-2}$ and $G$ is obtained from the complete balanced bipartite graph with parts $A$ and $B$ of size $\frac{n}{2}$ by adding an $\frac{n}{2(\ell-1)}$-regular graph in each part. Similarly as for triangles, to cover the $\frac{n^2}{4}$ edges between $A$ and $B$, one requires at least $\frac{n^2}{4(\ell-1)}$ copies of $C_\ell$, each using at least one of the $\frac{n^2}{4(\ell-1)}$ edges within $A$ or $B$. Thus, no graph obtained from a $G\in \G_\ell$ by removing edges within $A$ or $B$ can have a $C_\ell$-decomposition. In particular, when $n=2k\ell(\ell-1)$ for some odd natural number $k$ and a perfect matching is removed from each of $A$ and $B$, it is not too difficult to verify that the resulting graph is $C_\ell$-divisible, thus showing that minimum degree $(\frac{1}{2}+\frac{1}{2(\ell-1)})n$ is required to guarantee a $C_\ell$-decomposition of a $C_\ell$-divisible graph.
A natural conjecture states that this should suffice, and the asymptotic version is thus the following.

\begin{conj}\label{conj:odd}
    For all odd $\ell\geq 3$, we have $\delta_{C_\ell}=\frac{1}{2}+\frac{1}{2(\ell-1)}$.
\end{conj}

Progress towards this conjecture was first obtained by Barber, K\"uhn, Lo, and Osthus \cite{BarKuhLoOst}, who showed that $\delta_{C_\ell}\leq \delta_{K_3}$
for any odd $\ell$. This was then significantly improved by Joos and K\"uhn \cite{JooKuh}, who showed that, for any $\varepsilon$, there exists $\ell_0$ such that $\delta_{C_\ell}\leq \frac{1}{2}+\varepsilon$ for any odd $\ell\geq \ell_0$.
Finally, Bryant, Dukes, Horsley, Maenhaut, and Montgomery \cite{BryDukHorMaeMon} gave an explicit formula for this $\varepsilon$-term by showing that $\delta_{C_\ell} \leq \frac{1}{2}+\frac{1}{2(\ell-1)-2}$ for any odd $\ell\geq 5$. As previously discussed, the case $\ell=3$ was resolved recently independently of this paper (\cite{DelPos2026}, see \cref{conj:NW}).
We resolve \cref{conj:odd} for any odd $\ell \geq 73$.
 
\begin{theorem}\label{T:main}
    For all odd $\ell\geq 73$, we have $\delta_{C_{\ell}} = \frac{1}{2}+\frac{1}{2(\ell-1)}$.
\end{theorem}

Note that the bound of $73$ is an artefact of our proof, and we have no reason to believe that the conjecture is not also true for all smaller odd values of $\ell$ (especially since it is now known to be true for $\ell=3$). We only make moderate attempts at optimising this bound. 
Indeed, on the one hand, it is clear that new ideas would be needed to obtain the case $\ell=5$.
On the other hand, note that if one were only interested in proving the existence of an $\ell_0$ such that \cref{conj:odd} holds for all $\ell\geq \ell_0$, then our proof could be simplified and shortened. We thus aim for a balance between technicality and obtaining an explicit sensible value of $\ell$.
This is discussed in more detail in \cref{sec:conclusion}.

The paper is organised as follows. The proof of \cref{T:main} is divided into a non-extremal (\cref{sec:nonextremal}) and an extremal case (\cref{sec:extremal}); see \cref{sec:sketch} for a detailed proof overview and \cref{sec:conclusion} for a derivation of \cref{T:main}. Our notation and preliminary results are collected in \cref{sec:preliminaries}, and we finish with some concluding remarks in \cref{sec:conclusion}.

\section{Proof overview}\label{sec:sketch}

We now give an overview of the proof of our main result, \cref{T:main}. Let $\ell\geq 73$ and $n$ be sufficiently large. Let $G$ be an $n$-vertex graph with minimum degree $\delta(G)\geq (\frac{1}{2}+\frac{1}{2(\ell-1)})n$. Recall that by \cref{T:reduction}, it suffices to show that $G$ has an approximate $C_\ell$-decomposition.

\subsection{General strategy}

In \cite{BryDukHorMaeMon}, it was shown that in order to find an approximate $C_\ell$-decomposition of $G$, it in fact suffices to find a fractional decomposition into smaller weighted structures called `condensations'. We will decompose into two different condensations of $C_\ell$. Both are formed by beginning with a `base' triangle with weight 1 on each edge and adding an edge of weight $\ell-3$. The added edge may be incident to two vertices of the triangle (in which case its weight is added to form an edge of weight $\ell-2$ in the triangle) or incident to one vertex of the triangle (in which case it becomes a pendant edge of weight $\ell-3$). See \cref{sec:cond} for details.

Our overall strategy is therefore as follows. We first find a fractional $K_3$-packing $\cD$ of $G$, which will form the base triangles of our condensations.
Let $L$ be the `leftover' weighted subgraph of $G$. To each triangle in $\cD$, we wish
to add (up to) a certain weight of incident edges to form our desired condensations.
That is, it suffices to show that the leftover weights in $L$ can be suitably distributed among the triangles in $\cD$. 
In practice, we will see in \cref{L:hallCondition} that we need to verify that the leftover weight is not too large in total and also well distributed across the entire vertex set. More precisely, we need to check a condition in the style of Hall's marriage theorem: for any $X\subseteq V(G)$, the total leftover weight in $L[X]$ is not too large compared to the combined weights of triangles in $\cD$ which are incident to $X$.

\subsection{Finding a suitable fractional triangle packing}

The construction of the fractional $K_3$-packing $\cD$  varies
depending on the structure of $G$.%
\footnote{Note that we require our packing to obey local properties, and so a result that only guarantees the existence of a packing of a certain size, like the recent result in \cite{DelPos2026}, is not sufficient for our purposes.}

A na\"ive way to create $\D$ is to have each edge offer an equal share of its weight to each of the triangles that contains it, to set the coefficient of each triangle to be the least of the three shares it is offered, and then increase coefficients arbitrarily until the packing is maximal. Intuitively, this approach should produce a suitable packing provided that the distribution of the numbers of triangles using each edge of $G$ is not too extreme.%
\footnote{In fact, this was the strategy used in \cite{BryDukHorMaeMon} to show that $\delta_{C_{\ell}} \leq  \frac{1}{2}+\frac{1}{2(\ell-1)-2}$. This minimum degree condition essentially corresponds to the threshold where $\cD$ immediately translates into a fractional decomposition into weighted triangles with two edges of weight $1$ and one edge of weight $\ell-2$. The core of the work in \cite{BryDukHorMaeMon} was to show that a fractional decomposition into condensations gives rise to a $C_\ell$-decomposition.}
In fact, in \cref{sec:nonextremal}, we show that this na\"ive approach produces a packing which obeys the required Hall-type condition  whenever $G$ is reasonably far from the graphs in $\G_\ell$, in the sense that it has significantly more edges than the graphs in $\G_\ell$ or has no cut of size close to $\frac{n^2}{4}$. 

We may therefore assume that $G$ is \emph{near $\ell$-extremal}, that is, $G$ does not have many more edges than the graphs in $\G_\ell$ and also has an edge cut of size close to $\frac{n^2}{4}$.
To gain intuition, first consider the case where $G \in \G_\ell$. Let $\{A,B\}$ be the equipartition of $V(G)$ which induces a complete bipartite graph. 
Recall from the introduction that any $C_\ell$-decomposition of a graph in $G_\ell$ has a very specific structure. Similarly, in a fractional decomposition of $G$ into our two desired condensations, each copy of a condensation
must have at least $\frac{1}{\ell}$ of its weight lying on \emph{side edges} within $A$ or $B$.
But exactly $\frac{1}{\ell}$ of the edges of $G$ are side edges, so to obtain the desired decomposition, each copy of a condensation
must have must have \emph{exactly} $\frac{1}{\ell}$ of its weight lying on side edges.
It is not too difficult to see that this can only be achieved if each `base' triangle uses only one side edge and the added edge is not a side edge. In other words, we require that
$\D$ uses each side edge completely and that each of its triangles must contain only one side edge. Let $\D$ be the fractional packing in which each such triangle is taken with coefficient $\frac{2}{n}$. Since each of the side edges belongs to exactly $\frac{n}{2}$ such triangles, they are indeed fully covered.
On the other hand, each \emph{cross edge} between $A$ and $B$ belongs to precisely $\frac{n}{\ell-1}$ triangles in $G$ and so has $1-\frac{2}{n}\cdot \frac{n}{\ell-1}=\frac{\ell-3}{\ell-1}$ weight in the leftover $L$. It is not too hard to show that this $\D$ will obey the required Hall-type condition.

When $G$ is only near $\ell$-extremal, there is still a partition $\{A,B\}$ of $V(G)$ with $|A|\approx |B|$ and all but a small set $U$ of \emph{exceptional vertices} have degrees similar to those in graphs in $\G_\ell$ (degree roughly $\frac{n}{2(\ell-1)}$ to their own part and degree roughly $\frac{n}{2}$ to the other part). In this case our approach to creating $\D$ is similar to our na\"\i ve one, except that we entirely ignore triangles within $A$ and $B$. We show that there is then only a little leftover weight within $A$ and $B$, and a typical cross edge can only have slightly more than $\frac{\ell-3}{\ell-1}$ leftover weight. Intuitively, it should be that $\cD$ is overall `better' than the fractional $K_3$-packing of a graph in $\G_\ell$ defined above. Indeed, in $G$, there are still at most $\frac{n^2}{4}$ edges between $A$ and $B$, so any deviation from graphs in $\G_\ell$ (more edges in total and/or smaller maximum cut) gives rise to a larger proportion of side edges versus cross edges, and so, a larger proportion of edges with low leftover weight versus edges with high leftover weight. 
We show in \cref{sec:extremal} that this intuition holds true and that the Hall-type condition holds for $\D$.

Our method of showing that the desired Hall-type condition holds varies depending on the situation: 
when $G$ is not near $\ell$-extremal (\cref{lm:nonextremal}), we utilise the fact that $\cD$ must fully cover relatively many edges; when $G$ is near $\ell$-extremal and $X$ is not too large (\cref{lm:extremalsmall}), we analyse the weight on cross edges used by $\cD$; when $G$ is near $\ell$-extremal and $X=V(G)$ (\cref{lm:extremalfull}), we analyse the weight on side edges used by $\cD$; and when $G$ is near $\ell$-extremal and $X\neq V(G)$ but is very large (\cref{lm:extremalXbar}), we make use of an additional technical property of our packing which allows us to control the influence of exceptional vertices. We provide more specific details before the proof of each lemma.

\section{Preliminaries}\label{sec:preliminaries}

\subsection{Notation}

The following notational conventions will be used throughout the paper.

\subsubsection{Standard notation}
Given a positive 
integer $k$, we let $[k]\coloneqq \{1,\dots, k\}$. We denote by $\mathbb{R}^+$ the set of positive reals and by $\mathbb{R}^{\geq 0}$ the set of non-negative reals. 

A graph consists of a set $V$ of \emph{vertices} and a set $E\subseteq \binom{V}{2}$ of \emph{edges}.
Let $G$ be a graph. We denote by $V(G)$ and $E(G)$ its vertex and edge sets, and define $e(G)\coloneqq |E(G)|$. The \emph{cycle} on $\ell$ vertices and edges is denoted by $C_\ell$ and the \emph{complete graph} on $n$ vertices is denoted by $K_n$.

Given a vertex $v\in V(G)$, we denote by $N_G(v)$ its \emph{neighbourhood} and define its \emph{degree} as $d_G(v)\coloneqq |N_G(v)|$. The \emph{minimum degree} of $G$ is $\delta(G)\coloneqq \min_{v\in V(G)}d_G(v)$ and the \emph{maximum degree} of $G$ is $\Delta(G)\coloneqq \max_{v\in V(G)}d_G(v)$. Given disjoint $A,B\subseteq V(G)$, we denote by $E_G(A,B)$ the set of edges of $G$ with one endpoint in $A$ and one endpoint in $B$, we let $e_G(A,B)\coloneqq |E_G(A,B)|$, and we denote by $G[A,B]$ the bipartite graph on vertex classes $A$ and $B$ and edge set $E_G(A,B)$. 

Let $X\subseteq V(G)$. We denote by $G[X]$ the subgraph of $G$ \emph{induced} by $X$, that is, the graph on vertex set $X$ where two vertices $u,v\in X$ are adjacent in $G[X]$ if and only if they are adjacent in $G$.
We say $X$ is \emph{independent} if $e(G[X])=0$. 

The \emph{complement} of $G$ is the graph $\overline{G}$ on $V(G)$ where two vertices $u,v\in V(G)$ are adjacent in $\overline{G}$ if and only if they are not adjacent in $G$. For any set $X$ of vertices of a graph $G$, we use $\overline{X}$ to denote $V(G) \setminus X$ (the graph $G$ will be clear from the context and so this will be well defined).

\subsubsection{Flow networks}
A directed graph consists of a set $V$ of \emph{vertices} and a set $E\subseteq \{(u,v) \in V \times V: u\neq v\}$ of \emph{arcs}. 
Let $G$ be a directed graph. We denote by $V(G)$ and $E(G)$ its vertex and arc sets, and given $A,B\subseteq V(G)$, we define $E_G(A,B)\coloneqq \{(u,v)\in E(G): u\in A, v\in B\}$.

A \emph{network} $N$ is a directed graph, together with a \emph{source} $s_1\in V(N)$, a \emph{sink} $s_2\in V(N)$, and a \emph{capacity} function $c: E(N) \rightarrow \mathbb{R}^+ \cup \{\infty\}$.
A \emph{flow} in a network $N$ is then a function $f:E(N) \rightarrow \mathbb{R}^+$ such that $f(e)\leq c(e)$ for each $e\in E(N)$ and $\sum_{u: (u,v)\in E(N)}f((u,v))=\sum_{w: (v,w)\in E(N)}f((v,w))$ for each $v\in V(N)\setminus \{s_1,s_2\}$.
The \emph{value} of a flow $f$ is equal to $\sum_{v: (s_1,v)\in E(N)}f((s_1,v))=\sum_{v: (v,s_2)\in E(N)}f((v,s_2))$.
A \emph{cut} in a network $N$ is a partition $(S_1,S_2)$ of $V(N)$ such that $s_1\in S_1$ and $s_2\in S_2$. The \emph{capacity} of a cut $(S_1,S_2)$ is equal to $\sum_{(u,v)\in E_N(S_1,S_2)}c((u,v))$.

\subsubsection{Weighted graphs}
A \emph{weighted graph} $W$ consists of an \emph{underlying graph} $H$, together with an assignment $w_W$ of positive weights to the edges of $H$. We use $V(W)$, $E(W)$, and $e(W)$ to denote $V(H)$, $E(H)$, and $e(H)$. We let $\|W\|\coloneqq \sum_{e\in E(W)}w_W(e)$ denote the sum of the weights on the edges of $W$ and, for a set $\mathcal{W}$ of weighted graphs, we define $\|\mathcal{W}\| \coloneqq \sum_{W \in \mathcal{W}}\|W\|$. For convenience we let $w_W(xy)=0$ for any $xy \notin E(H)$. We often consider graphs as weighted graphs in which each edge has weight 1.

Let $W$ and $W'$ be weighted graphs. We define $W \uplus W'$ to be the weighted graph such that $V(W \uplus W')=V(W) \cup V(W')$, $E(W \uplus W')=E(W) \cup E(W')$, and $w_{W \uplus W'}(e)=w_W(e)+w_{W'}(e)$ for each $e \in E(W \uplus W')$. We say $W'$ is a \emph{weighted subgraph} of $W$ if the underlying graph of $W'$ is a subgraph of the underlying graph of $W$ and $w_{W'}(e) \leq w_W(e)$ for all $e \in E(W')$. If $W'$ is a weighted subgraph of $W$, we define $W-W'$ to be the weighted subgraph of $W$ such that, for each $e \in E(W)$, we have $e \notin E(W-W')$ if $w_{W'}(e) = w_W(e)$ and $w_{W-W'}(e) = w_W(e)-w_{W'}(e)$ otherwise.

\subsubsection{Fractional decompositions and packings}
For weighted graphs $W_0$ and $W$, we say $W_0$ is a \emph{scaled copy} of $W$ if there is a real number $\alpha$ and an isomorphism $f$ from the underlying graph of $W$ to the underlying graph of $W_0$ such that $w_{W_0}(e)=\alpha\, w_W(f(e))$ for all $e \in E(W)$. This (unique) number $\alpha$ is called the \emph{coefficient}.

For a set $\mathcal{W}$ of weighted graphs, a \emph{fractional $\mathcal{W}$-decomposition} of a weighted graph $G$ is a finite collection $\mathcal{D}=\{D_1,\ldots,D_t\}$ of weighted subgraphs of $G$, each a scaled copy of a weighted graph in $\mathcal{W}$, such that $G=D_1 \uplus \cdots \uplus D_t$. A \emph{fractional $\mathcal{W}$-packing} of $G$ is a fractional $\mathcal{W}$-decomposition $\D$ of a weighted subgraph $D$ of $G$, and its \emph{leftover} is the weighted graph $G-D$. Given $e\in E(D)$, we write $w_{\D}(e)$ for $w_D(e)$. For a fractional packing $\D$ of a graph $G$ and a subset $X$ of $V(G)$, we define $\D_X\coloneqq \{D_i \in \D:V(D_i) \cap X \neq \emptyset\}$ to be the fractional packing induced from $\cD$ by only considering the weighted subgraphs incident to $X$. When $\mathcal{W}$ contains a single weighted graph $W$, we replace $\mathcal{W}$ by $W$ in the above notation.

\subsubsection{Triangles}
The \emph{triangle} is the graph $K_3=C_3$.

Let $G$ be a graph. We denote by $\T(G)$ the set of all triangles in $G$. Given $e\in E(G)$, we let $\T_G(e)$ be the set of triangles in $\T(G)$ which contain $e$ and we let $t_G(e)\coloneqq |\T_G(e)|$ denote the number of triangles in $G$ that $e$ is in.

By a slight abuse of notation, we will sometime define a fractional $K_3$-packing of $G$ by a function $w_{\D}:\T(G) \rightarrow \mathbb{R}^{\geq 0}$. 
Provided that for each $e\in E(G)$, we have $\sum_{T\in \T_G(e)}w_{\D}(T)\leq 1$, such a function defines the fractional packing $\D\coloneqq \{W_T:T \in \T(G), w_{\D}(T)>0\}$ where $W_T$ is the scaled copy of $T$ with coefficient $w_{\D}(T)$.

\subsubsection{Near extremal graphs}

Let $\ell\geq 5$ be odd. An $n$-vertex graph $G$ is called \emph{near $\ell$-extremal} if 
\begin{enumerate}[label={(E\arabic*)}]
    \item \label{extremalDeg}
$\delta(G)\geq \frac{1}{2}(1+\frac{1}{\ell-1})n$;
    \item \label{extremalSize}
$e(G) \leq \frac{1}{4}(1+\frac{1}{\ell-1} +\frac{2}{(\ell-3)(\ell-1)})n^2$; and
    \item \label{extremalCut}
there is a partition $\{A,B\}$ of $V(G)$ with $e_G(A,B)\geq (\frac{1}{4}-\frac{1}{(\ell-3)(\ell-1)})n^2$.
\end{enumerate}

\subsection{Mantel's theorem}

We will need some standard results about triangle-free graphs.

\begin{lemma}[Mantel \cite{mantel1907solution}]\label{lm:Mantel}
    Any $n$-vertex triangle-free graph has at most $\frac{n^2}{4}$ edges.
\end{lemma}

The following well-known result gives a stability version of Mantel's theorem. We give the short proof for the sake of completeness. 

\begin{lemma}\label{L:MantelStability}
Let $H$ be an $n$-vertex triangle-free graph. Then $H$ has a bipartite subgraph with at least $2e(H)-\frac{n^2}{4}$ edges.
\end{lemma}

\begin{proof}
Let $z$ be a vertex of maximum degree $\Delta(H)$ in $H$ and let $S \coloneqq V(H) \setminus N_H(z)$. Since $N_H(z)$ is independent, the edges in $E(H) \setminus E(H[S])$ form a bipartite subgraph of $H$. Now
\[|E(H)|+|E(H[S])| = \sum_{x \in S}d_H(x) \leq |S|\Delta(H) =  |S||\overline{S}| \leq \mfrac{n^2}{4}\]
and hence $|E(H[S])| \leq \frac{n^2}{4}-|E(H)|$. Thus $|E(H)|-|E(H[S])| \geq 2|E(H)|-\frac{n^2}{4}$.
\end{proof}

\subsection{Condensations}\label{sec:cond}

We will now define condensations formally and state the key result of \cite{BryDukHorMaeMon}.

Let $F$ be a graph. Given a partition $\P=\{V_1, \dots , V_k\}$ of $V(F)$ into independent sets,
the \emph{condensation} of $F$ with respect $\P$ 
is the weighted graph $W$ on vertex set $[k]$ where the edge $ij$ is absent if $e_F(V_i,V_j)=0$, and the edge $ij$ is present with weight $e_F(V_i,V_j)$ otherwise. We say a weighted graph $W$ on vertex set $[k]$ is a \emph{condensation} of $F$ if there exists partition $\mathcal{P}$ of $V(F)$ such that $W$ is the condensation of $F$ with respect to $\mathcal{P}$.

Bryant, Dukes, Horsley, Maenhaut, and Montgomery \cite{BryDukHorMaeMon} showed that to find an approximate $F$-decomposition, it suffices to find a fractional decomposition into condensations of $F$.

\begin{theorem}[Bryant, Dukes, Horsley, Maenhaut, and Montgomery \cite{BryDukHorMaeMon}]\label{T:condensations}
Let $F$ be a graph, let $\mathcal{W}$ be a set of condensations of $F$, and let $\eta>0$ be a real number. There is an integer $n_0 \coloneqq n_0(F,\eta)$ such that the following holds for all $n\geq n_0$. Let $G$ be an $n$-vertex graph. If $G$ has a fractional $\mathcal{W}$-decomposition, then it also has an $\eta$-approximate $F$-decomposition.
\end{theorem}

As discussed in \cref{sec:sketch}, we will restrict ourselves to the following two specific condensations. Let $C^{\stri}_{\ell}$ denote the weighted triangle with two edges of weight $1$ and one edge of weight $\ell-2$ and let $C^{\splus}_{\ell}$ be the weighted graph whose underlying graph is a triangle with a pendant edge, where the pendant edge has weight $\ell-3$, and all edges of the triangle have weight $1$.

\begin{lemma}\label{lm:cond}
    Let $\ell\geq 5$ be odd. Then $C^{\stri}_{\ell}$ and $C^{\splus}_{\ell}$ are condensations of $C_\ell$.
\end{lemma}

\begin{proof}
Let $V(C_\ell)=\{v_i: i\in [\ell]\}$ and $E(C_\ell)=\{v_iv_{i+1}: i\in [\ell]\}$ (where $v_{\ell+1}\coloneqq v_1$). Then, 
the partition $\{\{v_1\},\{v_2,v_4,\ldots,v_{\ell-1}\},\{v_3,v_5,\ldots,v_{\ell}\}\}$ witnesses that $C^{\stri}_{\ell}$ is a condensation of $C_\ell$, while the partition $\{\{v_1\},\{v_2\},\{v_3,v_5,\ldots,v_{\ell}\},\{v_4,v_6,\ldots,v_{\ell-1}\}\}$ witnesses that $C^{\splus}_{\ell}$ is a condensation of $C_\ell$.
\end{proof}

\subsection{A simple fractional triangle packing}

As discussed in \cref{sec:sketch}, to find a fractional decomposition into our two condensations $C^{\stri}_{\ell}$ and $C^{\splus}_{\ell}$, we will start by constructing a fractional $K_3$-packing. In this section, we describe a basic construction of such a packing, which will be used in the non-extremal case.

Recall that $t_G(e)$ denotes the number of triangles of $G$ that contain the edge $e$.
It was observed in \cite{BryDukHorMaeMon} that, if two edges $e,e'$ are adjacent, then $t_G(e)$ and $t_G(e')$ cannot differ by too much. We include the short proof for completeness.

\begin{lemma}[\cite{BryDukHorMaeMon}]\label{lm:ratio}
Let $\ell \geq 3$ be an integer and let $G$ be an $n$-vertex graph with minimum degree at least $(\frac{1}{2}+\frac{1}{2(\ell-1)})n$. For any adjacent edges 
$e,e'\in E(G)$, we have $t_G(e) \leq \frac{\ell}{2}t_G(e')$.
\end{lemma}

\begin{proof}
Note that $t_G(e') \geq 2\delta(G)-n \geq \frac{n}{\ell-1}$, so we may assume that $t_G(e)>\frac{\ell n}{2(\ell-1)}$. 
Since any vertex 
is nonadjacent to at most $\frac{(\ell-2)n}{2(\ell-1)}$ vertices in $G$, we have
\[\mfrac{t_G(e)}{t_G(e')} \leq \mfrac{t_G(e)}{t_G(e)-\frac{(\ell-2)n}{2(\ell-1)}} < \mfrac{\ell}{2}\]
where the last inequality holds since $t_G(e) > \frac{\ell n}{2(\ell-1)}$ by assumption and the middle expression is decreasing in $t_G(e)$.
\end{proof}

This immediately gives a fractional $K_3$-packing where no edge has very small weight.

\begin{cor}\label{L:decentPartialDecomp}
Let $G$ be an $n$-vertex graph with minimum degree at least $(\frac{1}{2}+\frac{1}{2(\ell-1)})n$. There is a fractional $K_3$-packing $\cD$ of $G$ such that $w_\cD(e)\geq \frac{2}{\ell}$ for each $e\in E(G)$.
\end{cor}

\begin{proof}
Let $\D$ be the fractional $K_3$-packing of $G$ given by $w_{\D}(T) \coloneqq \min\{\frac{1}{t_G(e)}:e \in E(T)\}$ for each $T \in \T(G)$. Let $e \in E(G)$. Since by \cref{lm:ratio}, $t_G(e') \leq \frac{\ell}{2}t_G(e)$ for each $e'\in E(G)$ which is adjacent to $e$, we have $w_{\D}(T) \geq \frac{2}{\ell t_G(e)}$ for each $T \in \T_G(e)$. Thus $w_{\D}(e) \geq \frac{2}{\ell}$.
\end{proof}

\subsection{Absorbing the leftover}

Recall from \cref{sec:sketch} that after constructing a fractional $K_3$-packing $\cD$, we need to distribute the leftover to the triangles in $\cD$ to obtain copies of our condensations $C^{\stri}_{\ell}$ and $C^{\splus}_{\ell}$.

The following lemma states that if $T$ is a triangle in $\cD$, then we can assign it up to a certain amount of leftover weight to form copies of $C^{\stri}_{\ell}$ and $C^{\splus}_{\ell}$. 

\begin{lemma}\label{L:hairyTriDec}
Let $T$ be a scaled copy of $K_3$ and let $L$ be a weighted graph such that $\|L\| \leq \frac{\ell-3}{3}\|T\|$ and each edge of $L$ has at least one endpoint in $V(T)$. Then $T \uplus L$ has a fractional $\{C^{\stri}_{\ell},C^{\splus}_{\ell}\}$-decomposition.
\end{lemma}

\begin{proof}
Since $T$ is a scaled triangle, recall that we have $w_T(e)=\frac{1}{3}\|T\|$ and so, by assumption, $w_T(e)\geq \frac{1}{\ell-3}\|L\|$ for all $e\in E(T)$.
For each $e\in E(L)$ that has exactly one endpoint in $V(T)$, let $W_e$ be the weighted copy of $C^{\splus}_{\ell}$ in which $e$ has weight $w_L(e)$ and each edge of $T$ has weight $\frac{1}{\ell-3}w_L(e)$. For each $e\in E(L)$ that has both endpoints in $V(T)$, let $W_e$ be the weighted copy of $C^{\stri}_{\ell}$ in which $e$ has weight $\frac{\ell-2}{\ell-3}w_L(e)$ and each other edge of $T$ has weight $\frac{1}{\ell-3}w_L(e)$. Then $\sum_{e \in E(L)}w_{W_e}(e')=w_L(e')$ for each $e' \in E(L) \setminus E(T)$ and, for each $e' \in E(T)$, 
\[\sum_{e \in E(L)}w_{W_e}(e')=\mfrac{\ell-2}{\ell-3}w_L(e')+\mfrac{1}{\ell-3}\bigl(\|L\|-w_L(e')\bigr)=\mfrac{1}{\ell-3}\|L\|+w_L(e')\leq w_T(e')+w_L(e').\]
So $\{W_e:e \in E(L)\}$ is a fractional $\{C^{\stri}_{\ell},C^{\splus}_{\ell}\}$-packing of $T \uplus L$. The leftover of this packing is empty if $\|L\| = \frac{\ell-3}{3}\|T\|$ and otherwise its edges form a scaled copy of $K_3$ (with coefficient $\frac{1}{3}\|T\|-\frac{1}{\ell-3}\|L\|$), which clearly has a fractional $C^{\stri}_{\ell}$-decomposition. So we are done in either case.
\end{proof}

Using the max-flow min-cut theorem, we derive a Hall-type condition which guarantees that the leftover of a fractional $K_3$-packing can be assigned in such a way that  \cref{L:hairyTriDec} can be applied to each of the resulting weighted graphs, and thus the desired fractional $\{C^{\stri}_{\ell},C^{\splus}_{\ell}\}$-decomposition of $G$ can be obtained.

\begin{theorem}[Max-flow min-cut \cite{ForFul}]\label{thm:MFMC}
    The maximum value of a flow in a network $N$ equals the minimum capacity of cut in $N$.
\end{theorem}

Recall that if $\cD$ is a fractional packing and $X$ is a vertex set, then $\cD_X$ is the fractional packing comprising all the elements of $\cD$ that are incident to $X$.

\begin{lemma}\label{L:hallCondition}
Let $G$ be a graph, let $\D$ be a fractional $K_3$-packing of $G$, and let $L$ be the leftover of $\D$.
Suppose that for each $X \subseteq V(G)$, we have
\begin{equation}\label{E:hallCondition}
\mfrac{\ell-3}{3}\|\D_X\|-\|L[X]\| \geq 0.
\end{equation}
Then $G$ has a fractional $\{C^{\stri}_{\ell},C^{\splus}_{\ell}\}$-decomposition.
\end{lemma}

\begin{proof}
We set up an auxiliary flow network $N$ as follows. The vertex set of $N$ is $\{s_1,s_2\} \cup E(G) \cup \D$ where $s_1 \notin E(G) \cup \D$ is the source and $s_2 \notin E(G) \cup \D$ is the sink. For each $e \in E(G)$ there is an arc of capacity $w_L(e)$ from $s_1$ to $e$, and there is an arc of infinite capacity from $e$ to each $T \in \D$ such that $e$ is incident with at least one vertex of $T$. For each $T \in \D$, there is an arc of capacity $\frac{\ell-3}{3}\|T\|$ from $T$ to $s_2$. Denote by $c$ the capacity function in $N$. Note that sum of the capacities of the arcs leaving $s_1$ is $\|L\|$. We claim that every cut of $N$ has capacity at least $\|L\|$ across it.

Let $(S_1,S_2)$ be a cut in $N$ where $s_1 \in S_1$ and $s_2 \in S_2$. Let $X$ be the set of vertices of $G$ incident to at least one edge in $E(G) \cap S_1$. If there exists $T\in \cD_X\cap S_2$, then there is an arc of infinite capacity across $(S_1,S_2)$. We may therefore assume that $\D_X \subseteq S_1$.
Then the capacity across $(S_1,S_2)$ is
\begin{align*}
    \sum_{e\in E(G)\cap S_2}c((s_1,e))+\sum_{T\in \cD\cap S_1}c((T,s_2))&=
    \sum_{e \in E(G) \cap S_2}w_L(e)+\sum_{T \in \D \cap S_1}\mfrac{\ell-3}{3}\|T\|\\
    &\geq \left(\|L\|-\sum_{e \in E(G) \cap S_1}w_L(e)\right)+\sum_{T \in \D_X}\mfrac{\ell-3}{3}\|T\|\\
    &\geq \|L\|-\|L[X]\| + \mfrac{\ell-3}{3}\|\D_X\| \geq \|L\|.
\end{align*}
Thus the network admits a flow of value $\|L\|$ by \cref{thm:MFMC}. 

Let $f$ be a flow of value $\|L\|$ and note we must have $f((s_1,e))=w_L(e)$ for each $e\in E(G)$.
For each $T \in \D$, let $E_T$ be the set of all $e \in E(G)$ for which $f((e,T))>0$,
and let $L_T$ be the weighted graph with edge set $E_T$ in which $w_{L_T}(e)\coloneqq f((e,T))$ 
for each $e \in E_T$. For each $e \in E(G)$, the preservation of flow in $N$ at vertex $e$ guarantees that $\sum_{T \in \D}w_{L_T}(e)=w_L(e)$. So $\{L_T:T \in \D\}$ is a fractional decomposition of $L$ and hence $\{T \uplus L_T:T \in \D\}$ is a fractional decomposition of $G$. For each $T \in \D$, the preservation of flow in $N$ at the vertex $T$ guarantees that $\|L_T\| \leq \frac{\ell-3}{3}\|T\|$ and thus $T \uplus L_T$ has a fractional $\{C^{\stri}_{\ell},C^{\splus}_{\ell}\}$-decomposition by \cref{L:hairyTriDec}. Thus $G$ has a fractional $\{C^{\stri}_{\ell},C^{\splus}_{\ell}\}$-decomposition.
\end{proof}

Note that \cref{E:hallCondition} holds trivially when $X=\emptyset$. In the case $X=V(G)$, using $\|L\|=e(G)-\|\D\|$, we have that \cref{E:hallCondition} is equivalent to $\|\cD\| \geq \frac{3}{\ell}e(G)$.

\section{The non-extremal case}\label{sec:nonextremal}

In this section we deal with host graphs $G$ that are not near $\ell$-extremal. Our strategy is as follows. We start with the fractional $K_3$-packing $\cD$ given by \cref{L:decentPartialDecomp} and arbitrarily increase coefficients so that each triangle in $G$ has an edge which is fully covered by $\cD$. Let $X\subseteq V(G)$ and suppose for a contradiction that \cref{E:hallCondition} does not hold. First, it is not too difficult to show that the number of edges incident to $X$ which are fully covered by $\cD$ must be small, for otherwise \cref{E:hallCondition} would immediately hold, a contradiction. We proceed to show that $|X|\approx n$ and $G$ must be highly structured. The intuition behind this is that, since $L$ is triangle-free by the maximality of $\cD$, \cref{L:MantelStability} implies that when $e(L)$ is large, $L$ consists of a large bipartite graph with a few extra edges.
This means $G$ contains a large bipartite subgraph and $e(G)$ is relatively small. But this cannot be the case as $G$ is not near $\ell$-extremal by assumption, a contradiction.

\begin{lemma}\label{lm:nonextremal}
Let $\ell \geq 5$, and let $G$ be an $n$-vertex graph with
minimum degree at least $(\frac{1}{2}+\frac{1}{2(\ell-1)})n$ that is not near $\ell$-extremal. Then $G$ has a fractional $\{C^{\stri}_{\ell},C^{\splus}_{\ell}\}$-decomposition.
\end{lemma}

\begin{proof}
Let $\cD$ be a maximal fractional $K_3$-packing of $G$ with leftover $L$. By \cref{L:decentPartialDecomp}, we may assume that $w_\cD(e)\geq \frac{2}{\ell}$ for all $e\in E(G)$.
Note that maximality implies that each triangle of $G$ has an edge $e$ such that $w_\cD(e)=1$. In other words, $L$ is triangle-free and so, by Mantel's theorem (\cref{lm:Mantel}), 
\begin{equation}\label{eq:L}
    e(L[S]) \leq \mfrac{1}{4}|S|^2
\end{equation} 
for any $S\subseteq V(G)$.

Assume for a contradiction that $G$ does not have a fractional $\{C^{\stri}_{\ell},C^{\splus}_{\ell}\}$-decomposition. Then, \cref{L:hallCondition} implies that there is $X\subseteq V(G)$ such that $\|L[X]\| > \frac{\ell-3}{3}\|\mathcal{D}_X\|$. 
It will be convenient to denote by $f_1$ and $f_2$ the number of edges $e\in E(G)$ incident to exactly one, respectively two, vertices of $X$ such that $w_\cD(e)=1$ (i.e.\ which are `fully covered' by $\cD$).
Then,
\begin{equation*}
    \mfrac{\ell-2}{\ell}e(L[X]) \geq \|L[X]\| > \mfrac{\ell-3}{3}\|\mathcal{D}_X\| \geq  \mfrac{\ell-3}{3}\left(\mfrac{2}{\ell}e(L[X])+\mfrac{2}{\ell}e_L(X,\overline{X})+f_1+f_2\right), 
\end{equation*}
and rearranging gives
\begin{equation}\label{eq:f}
    f_1+f_2 \leq \mfrac{1}{\ell-3}e(L[X])-\mfrac{2}{\ell}e_L(X,\overline{X})\stackrel{\text{\cref{eq:L}}}{\leq} \mfrac{1}{4(\ell-3)}|X|^2-\mfrac{2}{\ell}e_L(X,\overline{X}).
\end{equation}

For $i \in [3]$, 
let $\D_{X,i}$ be the fractional packing induced from $\cD$ by all triangles that contain exactly $i$ vertices in $X$ and observe that $\D_X=\D_{X,1} \cup \D_{X,2} \cup \D_{X,3}$.
Then,
\begin{align*}
    \mfrac{4}{3}\|\D_X\|+\mfrac{2}{3}\|\D_{X,3}\|&\geq \mfrac{2}{3}\|\D_{X,1}\|+\mfrac{4}{3}\|\D_{X,2}\|+2\|\D_{X,3}\|
    =\sum_{x \in X}\sum_{y\in N_G(x)}w_\cD(xy)\\
    &\geq \mfrac{2}{\ell}\left(\sum_{x \in X}d_G(x)-2f_2-f_1\right)+2f_2+f_1
    =\mfrac{2}{\ell}\sum_{x \in X}d_G(x)+\mfrac{2\ell-4}{\ell}f_2+\mfrac{\ell-2}{\ell}f_1.
\end{align*}
Rearranging and using that $\|\D_{X,3}\| \leq 3f_2$ (as $L$ is triangle-free), we obtain that
\begin{align*}
\|\D_X\|
&\geq \mfrac{3}{4}\left(\mfrac{2}{\ell}\sum_{x \in X}d_G(x)+\mfrac{2\ell-4}{\ell}f_2+\mfrac{\ell-2}{\ell}f_1\right)-\mfrac{1}{2}\|\D_{X,3}\| \\
&\geq \mfrac{3}{4}\left(\mfrac{2}{\ell}\sum_{x \in X}d_G(x)+\mfrac{2\ell-4}{\ell}f_2+\mfrac{\ell-2}{\ell}f_1\right)-\mfrac{3}{2}f_2 = \mfrac{3}{2\ell}\sum_{x \in X}d_G(x) +\mfrac{3\ell-6}{4\ell}f_1 -\mfrac{3}{\ell}f_2 \\
&\geq \mfrac{3}{4(\ell-1)}|X|n+\mfrac{3\ell-6}{4\ell}f_1-\mfrac{3}{\ell}f_2. 
\end{align*}
By assumption, this implies that
\begin{align}
0 < \|L[X]\| - \mfrac{\ell-3}{3}\|\D_X\| &\leq \mfrac{\ell-2}{\ell}e(L[X]) - \mfrac{\ell-3}{3}\left(\mfrac{3}{4(\ell-1)}|X|n+\mfrac{3\ell-6}{4\ell}f_1-\mfrac{3}{\ell}f_2\right) \notag\\
&\stackrel{\text{\cref{eq:L}}}{\leq} \mfrac{\ell-2}{4\ell}|X|^2 -\mfrac{\ell-3}{4(\ell-1)}|X|n-\mfrac{(\ell-3)(\ell-2)}{4\ell}f_1+\mfrac{\ell-3}{\ell}f_2 \notag\\
&\stackrel{\text{\cref{eq:f}}}{\leq}  \mfrac{\ell-2}{4\ell}|X|^2 -\mfrac{\ell-3}{4(\ell-1)}|X|n-\mfrac{(\ell-3)(\ell-2)}{4\ell}f_1 \notag\\
&\qquad\qquad+\mfrac{\ell-3}{\ell}\left(\mfrac{1}{4(\ell-3)}|X|^2-\mfrac{2}{\ell}e_L(X,\overline{X})-f_1\right) \notag\\
&= \mfrac{|X|}{4}\left(\mfrac{\ell-1}{\ell}|X|-\mfrac{\ell-3}{\ell-1}n\right)-\mfrac{(\ell-3)(\ell+2)}{4\ell}f_1
-\mfrac{2(\ell-3)}{\ell^2}e_L(X,\overline{X}). \label{E:diffBound}
\end{align}
We now proceed to show that this leads to a contradiction.

\begin{enumerate}[wide,label=\textbf{Case \arabic*:}]
    \item $0 \leq |X| \leq \frac{\ell(\ell-3)}{(\ell-1)^2}n =(1-\frac{\ell+1}{(\ell-1)^2})n$. Then the first term, and hence the entire expression, on the right side of \cref{E:diffBound} is nonpositive and we have a contradiction.
    
    \item $(1-\frac{\ell+1}{(\ell-1)^2})n < |X| < \frac{2\ell-7}{2\ell-5}n=(1-\frac{2}{2\ell-5})n$. Since $\ell^2-3\ell-2 \geq (\ell-1)(\ell-3)$ for $\ell \geq 5$, we have in this case that
    \begin{align*}e_L(X,\overline{X}) &> \sum_{x \in \overline{X}}(d_{G}(x)-|\overline{X}|)-f_1 
    \geq (n-|X|)\left(\mfrac{\ell n}{2(\ell-1)}-(n-|X|)\right)-f_1 \\
    &> \mfrac{2}{2\ell-5}\left(\mfrac{\ell }{2(\ell-1)}-\mfrac{\ell+1}{(\ell-1)^2}\right)n^2-f_1
    = \mfrac{\ell^2-3\ell-2}{(\ell-1)^2(2\ell-5)}n^2-f_1
    \geq \mfrac{\ell-3}{(\ell-1)(2\ell-5)}n^2-f_1.
    \end{align*}
    Thus, \cref{E:diffBound} implies that
    \begin{align*}
        0&<\mfrac{n}{4}\left(\mfrac{(\ell-1)(2\ell-7)}{\ell(2\ell-5)}n-\mfrac{\ell-3}{\ell-1}n\right)-\left(\mfrac{(\ell-3)(\ell+2)}{4\ell}-\mfrac{2(\ell-3)}{\ell^2}\right)f_1
        -\mfrac{2(\ell-3)^2}{\ell^2(\ell-1)(2\ell-5)}n^2\\
        &=\mfrac{-72 + 41 \ell - 7 \ell ^2}{4\ell^2(\ell-1)(2\ell-5)}n^2 -\mfrac{(\ell-3)(\ell^2+2\ell-8)}{4\ell^2}f_1<0,
    \end{align*}
    where the final inequality holds since $-72 + 41 \ell - 7 \ell ^2<0$ and $\ell-3>0$ and $\ell^2+2\ell-8>0$ for all $\ell\geq 5$. Once again, this is a contradiction.

    \item $(1-\frac{2}{2\ell-5})n \leq |X| \leq n$. Note that, on this domain, $\frac{1}{2}(n-|X|)^2+\frac{1}{4(\ell-3)}|X|^2$ is a convex quadratic in $|X|$ whose maximum is $\frac{n^2}{4(\ell-3)}$ (attained on both ends on the interval).
    Moreover, there are at most $\binom{n-|X|}{2}$ edges $e\in E(G[\overline{X}])$ such that $w_\cD(e)=1$. 
    Thus, we have
    \begin{align*}
        \mfrac{n^2}{4}&\stackrel{\text{\eqmakebox[case3]{\cref{eq:L}}}}{\geq} e(L)\geq e(G)-\left(\mbinom{n-|X|}{2}+f_1+f_2\right)
        \stackrel{\text{\cref{eq:f}}}{\geq}e(G)-\mfrac{1}{2}(n-|X|)^2-\mfrac{1}{4(\ell-3)}|X|^2
        \geq e(G)-\mfrac{n^2}{4(\ell-3)}\\
        &\stackrel{\text{\eqmakebox[case3]{}}}{\geq} \left(\mfrac{1}{2}+\mfrac{1}{2(\ell-1)}\right)\mfrac{n^2}{2}-\mfrac{n^2}{4(\ell-3)}\geq \left(\mfrac{1}{4}-\mfrac{1}{2(\ell-3)(\ell-1)}\right)n^2.
    \end{align*}
    Rearranging the first line, we obtain that 
    \[e(G) \leq \left(\mfrac{1}{4}+\mfrac{1}{4(\ell-3)}\right)n^2 = \left(\mfrac{1}{4}+\mfrac{1}{4(\ell-1)} +\mfrac{1}{2(\ell-3)(\ell-1)}\right)n^2.\]
    But the second line and \cref{L:MantelStability} imply that the underlying graph of $L$, and so $G$, has an edge cut of size at least
    \[\left(\mfrac{1}{4}-\mfrac{1}{(\ell-3)(\ell-1)}\right)n^2.\]
    That is, $G$ is near $\ell$-extremal, a contradiction.
\end{enumerate}
This completes the proof of the lemma.
\end{proof}

\section{The extremal case}\label{sec:extremal}

In this section, we prove the following.

\begin{lemma}\label{lm:extremal}
Let $\ell\geq 73$, and let $G$ be a near $\ell$-extremal graph. Then $G$ has a fractional $\{C^{\stri}_{\ell},C^{\splus}_{\ell}\}$-decomposition.
\end{lemma}

\subsection{Notation and preliminaries}

To avoid an excess of notation, we will adopt the following conventions throughout this section. Given a fractional $K_3$-packing, $L$ will always denote its leftover. The fractional packing in question will always be clear from the context, so $L$ will be well defined.

Let $\ell\geq 5$ be odd. 
Note that a graph may have several partitions $\{A,B\}$ of its vertex set which witness that it is near $\ell$-extremal. We will assume throughout the section that any near $\ell$-extremal graph $G$ comes with a fixed arbitrary $(A,B)\coloneqq (A(G),B(G))$ which maximises $e_G(A,B)$ and where $|A|\geq |B|$. 

Let $G$ be a near $\ell$-extremal graph.
We denote $m\coloneqq m(G)\coloneqq \frac{1}{2}(|A|-|B|)$, so that $|A|=\frac{n}{2}+m$ and $|B|=\frac{n}{2}-m$. The \emph{side edges} of $G$ are the edges of $\Gs\coloneqq G[A]\cup G[B]$ and the \emph{cross edges} of $G$ are the edges of $\Gc\coloneqq G[A,B]$.
We denote by $U\coloneqq U(G)$ the set of vertices $v\in V(G)$ such that $d_{\Gs}(v)\geq \frac{n}{20}$. The vertices in $U$ are to be understood as \emph{exceptional} vertices, whose side degree is significantly larger than the side degree in the $\ell$-extremal graphs in $\G_\ell$. The graph $G$ will always be clear from the context and so the above notation $A$, $B$, $m$, $U$ will be well defined. 

By the maximality of $e_G(A,B)$, we have that $d_{\Gc}(x) \geq d_{\Gs}(x)$ and hence
\begin{equation}\label{eq:cU}
    d_{\Gc}(x) \geq \mfrac{\delta(G)}{2}\geq \mfrac{\ell n}{4(\ell-1)}> \mfrac{n}{4}
\end{equation}
for any $x \in V(G)$.
Moreover, by the definition of $U$, we have that 
\begin{equation}\label{eq:cgood}
    d_{\Gc}(x) \geq \delta(G)-\mfrac{n}{20} > \mfrac{9n}{20}
\end{equation}
for any $x \in\overline{U}$. We will frequently use these facts. We also compute useful upper bounds on $m$ and $|U|$.

\begin{lemma}\label{L:mAndU}
Let $\ell \geq 73$ and let $G$ be a near $\ell$-extremal $n$-vertex graph.
    \begin{enumerate}
        \item $m\leq n((\ell-3)(\ell-1))^{-1/2} < \frac{n}{\ell-3} \leq \frac{n}{70}$.\label{mBound}
        \item $|U| \leq \frac{105n}{(\ell-3)(\ell-1)}$. \label{uBound}
        \item $|U| \leq \frac{3n}{200}$. \label{uBound2}
    \end{enumerate}
\end{lemma}

\begin{proof}
For \cref{mBound}, note that
\[\left(\mfrac{1}{4}-\mfrac{1}{(\ell-3)(\ell-1)}\right)n^2\stackrel{\text{\cref{extremalCut}}}{\leq}e_G(A,B)\leq |A||B|=\left(\mfrac{n}{2}+m\right)\left(\mfrac{n}{2}-m\right)=\mfrac{n^2}{4}-m^2,\]
and so rearranging gives $m\leq n((\ell-3)(\ell-1))^{-1/2} < \frac{n}{\ell-3} \leq \frac{n}{70}$ for $\ell\geq 73$, as desired.

We now prove \cref{uBound} and \cref{uBound2}.
First note that
$d_{\Gs}(x) \geq \delta(G)-|B| \geq \frac{n}{2(\ell-1)}+m$ for each $x \in A$, $d_{\Gs}(x) \geq \delta(G)-|A| \geq \frac{n}{2(\ell-1)}-m$ for each $x \in B$, and $d_{\Gs}(x) \geq \frac{n}{20}$ for each $x \in U$. Hence, 
\begin{align*}
\mfrac{n^2}{2(\ell-1)}+\mfrac{3n^2}{(\ell-3)(\ell-1)}
&\stackrel{\text{\eqmakebox[ubound]{\cref{extremalSize},\cref{extremalCut}}}}{\geq} 2\left(e(G)-e_G(A,B)\right)=\sum_{x \in V(G)}d_{\Gs}(x)\\
&\stackrel{\text{\eqmakebox[ubound]{}}}{\geq} \mfrac{n}{20}|U|+\left(\mfrac{n}{2}+m-|U|\right)\left(\mfrac{n}{2(\ell-1)}+m\right)+\left(\mfrac{n}{2}-m\right)\left(\mfrac{n}{2(\ell-1)}-m\right)\\
&\stackrel{\text{\eqmakebox[ubound]{}}}{=} |U|\left(\mfrac{n}{20}-\mfrac{n}{2(\ell-1)}-m\right)+\mfrac{n^2}{2(\ell-1)}+2m^2.    
\end{align*}
Rearranging, we obtain
\begin{equation}\label{E:uBound}
|U| \leq \frac{\frac{3n^2}{(\ell-3)(\ell-1)}- 2m^2}{\frac{n}{20}-\frac{n}{2(\ell-1)}-m} \stackrel{\text{\cref{mBound}}}{\leq} \frac{\frac{3n^2}{(\ell-3)(\ell-1)}}{\frac{n}{20}-\frac{n}{2(\ell-1)}-\frac{n}{\ell-3}} \stackrel{\text{$\ell \geq 73$}}{\leq} \mfrac{3024n}{29(\ell-3)(\ell-1)} < \mfrac{105n}{(\ell-3)(\ell-1)}.    
\end{equation}
This establishes the coarse bound \cref{uBound}. To obtain the more precise bound \cref{uBound2}, notice that from \cref{E:uBound} we have
\[|U| \leq \frac{\frac{3n^2}{(\ell-3)(\ell-1)}- 2m^2}{\frac{n}{20}-\frac{n}{2(\ell-1)}-m} \stackrel{\text{$\ell \geq 73$}}{\leq} \frac{\frac{n^2}{1680}- 2m^2}{\frac{31n}{720}-m}.\]
Routine analysis shows that, considered as a function of $m$ on the domain $[0,\frac{n}{70}]$, this last expression is concave and achieves its maximum value of
\[\mfrac{(217 - \sqrt{39529})n}{1260} < \mfrac{3n}{200}\]
when $m=\frac{1}{5040}(217 - \sqrt{39529})$. This establishes \cref{uBound2}.
\end{proof}

\subsection{The fractional packing}

The first step is to find a suitable fractional $K_3$-packing. As discussed in the proof overview, the general intuition is to construct $\cD$ na\"ively as in the non-extremal case, but ignoring side triangles. Proceeding in this way would give a fractional $K_3$-packing which behaves essentially like the optimal one for the extremal graphs in $\G_\ell$. However, we introduce an additional technicality which will enable us to better control the influence of the exceptional vertices in $U$. More precisely, we also ignore triangles which are fully within $U$ and cap the maximum weight a triangle may have in $\cD$ (this is \cref{lm:extremalD}\cref{D:Tsmall,D:cross}). The chosen cap value is large enough that it does not affect the triangles in $\overline{U}$ and we still obtain the desired properties: most cross edges carry some weight (this is \cref{lm:extremalD}\cref{D:crossmin}), while side edges are almost fully covered (this is \cref{lm:extremalD}\cref{D:side1,D:sideU}). However, we now ensure that we do not deviate too much from this behaviour, in the sense that cross edges in $\overline{U}$ cannot accumulate too much weight (this is \cref{lm:extremalD}\cref{D:crossmax}). This will allow us to verify \cref{E:hallCondition} when $X$ is large but not the entire vertex set (\cref{lm:extremalXbar}).

\begin{lemma}\label{lm:extremalD}
    Let $\ell \geq 73$ and let $G$ be a near $\ell$-extremal $n$-vertex graph. Then, there exists a fractional $K_3$-packing $\cD$ of $G$ such that the following properties hold.
    \begin{enumerate}
        \item For all $T\in \mathcal{T}(\Gs) \cup \mathcal{T}(G[U])$, we have $w_{\D}(T)=0$.\label{D:cross}
        \item For all $T\in \mathcal{D}$, we have $w_{\D}(T)\leq \frac{4}{n}$. \label{D:Tsmall}
        \item For all $e\in E(\Gc[\overline{U}])$, we have $w_{\D}(e)\leq \frac{2}{5}$.\label{D:crossmax}
        \item For all $e\in E(\Gc) \setminus E(\Gc[U])$, we have $w_{\D}(e)\geq \frac{2}{\ell}$.\label{D:crossmin}
        \item For all $e\in E(\Gs[\overline{U}])$, we have $w_{\D}(e)=1$.\label{D:side1}
        \item For all $e\in E(\Gs) \setminus E(\Gs[U])$, we have $w_{\D}(e)\geq \frac{2}{3}$.\label{D:sideU}
    \end{enumerate}
\end{lemma}

\begin{proof}
    For each $e\in E(G)$, let $t'_G(e)$ be the number of triangles in $\mathcal{T}(G) \setminus \mathcal{T}(\Gs)$ which contain $e$. Let $\cD_0$ be the fractional $K_3$-packing of $G$ where for all $T=xyz\in \mathcal{T}(G)$ we define
    \[w_{\cD_0}(T)\coloneqq \begin{cases}
        0 & \text{if }T\in \mathcal{T}(\Gs) \cup \mathcal{T}(G[U])\\
        \min\left\{\mfrac{1}{t'_G(xy)},\mfrac{1}{t'_G(yz)}, \mfrac{1}{t'_G(xz)},\mfrac{4}{n}\right\} & \text{if }T\in \mathcal{T}(G) \setminus (\mathcal{T}(\Gs) \cup \mathcal{T}(G[U])).
    \end{cases}\]
    Let $\cD$ be a fractional $K_3$-packing of $G$ that is maximal subject to the conditions that $w_{\D}(T)=0$ if $w_{\D_0}(T)=0$ and otherwise $w_{\cD_0}(T)\leq w_{\cD}(T)\leq \frac{4}{n}$.

    We claim that all the desired properties are satisfied. By construction, \cref{D:cross,D:Tsmall} hold. In proving \cref{D:crossmax,D:side1,D:sideU,D:crossmin}, we often use the facts that we have $t_G'(xy) \geq d_{\Gc}(x)+d_{\Gc}(y)-(\frac{n}{2}+m)$ for any $xy \in E(\Gs)$ and $t_G(yz) =t_G'(yz) \leq d_{\Gs}(y)+d_{\Gs}(z)$ for any $yz
    \in E(\Gc)$. 
    In particular, the latter implies that $t_G(e) \leq \frac{n}{10}$ for each $e \in E(\Gc[\overline{U}])$, so \cref{D:crossmax} follows immediately from \cref{D:Tsmall}. 

    For \cref{D:crossmin}, let $xy\in E(\Gc) \setminus E(\Gc[U])$. Since $t_G(xy)=t'_G(xy) \geq 2\delta(G)-n \geq \frac{n}{\ell-1}$, we have
    \begin{align*}
        w_{\mathcal{D}}(xy)&\geq \sum_{z\in N_G(x)\cap N_G(y)}\min\left\{\mfrac{1}{t'_G(xy)},\mfrac{1}{t'_G(xz)},\mfrac{1}{t'_G(yz)},\mfrac{4}{n}\right\} \\ &\geq \mfrac{t_G(xy)}{\max_{z\in N_G(x)\cap N_G(y)}\{t_G(xy),t_G(xz),t_G(yz),n/4\}} \stackrel{\text{\cref{lm:ratio}}}{\geq} \mfrac{t_G(xy)}{\max\{\ell t_G(xy)/2,n/4\}} \geq \mfrac{2}{\ell},
    \end{align*}
    as desired.

    For \cref{D:side1}, let $xy\in E(\Gs[\overline{U}])$ and suppose for a contradiction that $w_{\D}(xy)<1$. Then, for each $z \in (N_{\Gc}(x)\cap N_{\Gc}(y))\setminus U$, \cref{D:crossmax} implies that $\max\{w_{\D}(xy),w_{\D}(xz),w_{\D}(yz)\} <1$, and hence $w_{\D}(xyz)=\frac{4}{n}$ by definition of $\D$. Thus,
    \[w_{\cD}(xy)\geq \bigl|\bigl(N_{\Gc}(x)\cap N_{\Gc}(y)\bigr)\setminus U\bigr|\cdot\mfrac{4}{n} \stackrel{\text{\cref{eq:cgood}}}{>} \left(2\cdot\mfrac{9n}{20}-\mfrac{n}{2}-m-|U|\right)\cdot\mfrac{4}{n} \stackrel{\text{\cref{L:mAndU}\cref{mBound},\cref{uBound2}}}{>} 1,\]
    a contradiction.

    Finally, for \cref{D:sideU}, let $xu\in E(\Gs)\setminus E(\Gs[U])$ with $x\notin U$ and $u\in U$. First observe that
    \begin{equation}\label{E:mutualNeighbourhoodSize}
    |N_{\Gc}(x)\cap N_{\Gc}(u)| \stackrel{\text{\cref{eq:cgood}}}{\geq} \mfrac{9n}{20}+d_{\Gc}(u)-\mfrac{n}{2} - m \stackrel{\text{\cref{L:mAndU}\cref{mBound}}}{\geq} d_{\Gc}(u)-\mfrac{9n}{140} \stackrel{\text{\cref{eq:cU}}}{\geq} \mfrac{13n}{70}.
    \end{equation}
    Let $z\in N_{\Gc}(x)\cap N_{\Gc}(u)$. By \cref{eq:cU} and since $d_{\Gs}(x) < d_{\Gs}(u) \leq d_{\Gc}(u)$, we have
    \begin{equation*}\label{E:convexSum0}
    w_{\mathcal{D}}(xuz)\geq \min\left\{\mfrac{1}{d_{\Gc}(u)},\mfrac{1}{d_{\Gs}(x)+d_{\Gs}(z)},\mfrac{1}{d_{\Gs}(u)+d_{\Gs}(z)},\mfrac{4}{n}\right\}\geq \mfrac{1}{d_{\Gc}(u)+d_{\Gs}(z)}.
        \end{equation*}    
    Let $\mu_z\coloneqq d_{\overline{G}[A,B]}(z)$ denote the `missing' cross degree of $z$. Note that 
    \[d_{\Gs}(z)=d_G(z)-d_{\Gc}(z)\leq d_G(z)-\left(\mfrac{n}{2}-m-\mu_z\right) \leq \mfrac{n}{2(\ell-1)}+m+\mu_z+\left(d_G(z)-\mfrac{\ell n}{2(\ell-1)}\right).\]
    Thus, 
    \begin{equation}\label{E:convexSum}
        w_{\mathcal{D}}(xu)\geq \sum_{z\in N_{\Gc}(x)\cap N_{\Gc}(u)}\mfrac{1}{d_{\Gc}(u)+\frac{n}{2(\ell-1)}+m+\mu_z+d_G(z)-\frac{\ell n}{2(\ell-1)}}.
    \end{equation}    
    Each term of this sum is a convex function of $\mu_z+(d_G(z)-\frac{\ell n}{2(\ell-1)})$, and so, as 
    \begin{align*}
        \sum_{z\in N_{\Gc}(x)\cap N_{\Gc}(u)}\left(\mu_z+\left(d_G(z)-\mfrac{\ell n}{2(\ell-1)}\right)\right)&
        \stackrel{\text{\eqmakebox[missing]{}}}{\leq} \bigl(|A||B|-e_G(A,B)\bigr)+\left(2e(G)-\mfrac{\ell n^2}{2(\ell-1)}\right)\\
        &\stackrel{\text{\eqmakebox[missing]{\cref{extremalSize},\cref{extremalCut}}}}{\leq} \mfrac{2n^2}{(\ell-1)(\ell-3)},
    \end{align*}
   we have that
    \begin{align*}
        w_{\mathcal{D}}(xu)
        &\stackrel{\text{\eqmakebox[minweight]{}}}{\geq} \frac{|N_{\Gc}(x)\cap N_{\Gc}(u)|}{d_{\Gc}(u)+\frac{n}{2(\ell-1)}+m+\frac{2n^2}{(\ell-1)(\ell-3)|N_{\Gc}(x)\cap N_{\Gc}(u)|}}\\
        &\stackrel{\text{\eqmakebox[minweight]{\cref{E:mutualNeighbourhoodSize}, \cref{L:mAndU}\cref{mBound}}}}{\geq} \hspace{1cm} \frac{d_{\Gc}(u)-\frac{9n}{140}}{d_{\Gc}(u)+\frac{n}{144}+\frac{n}{70}+\frac{2n}{70\cdot72\cdot(13/70)}}\\
        &\stackrel{\text{\eqmakebox[minweight]{\cref{eq:cU}}}}{\geq} \mfrac{12168}{17911}>\mfrac{2}{3},
    \end{align*}
    as desired.
\end{proof}

\subsection{Distributing the leftover edges to the triangles}

Let $G$ be a near $\ell$-extremal graph. It now remains to show that given any fractional $K_3$-packing $\mathcal{D}$ of $G$ that obeys the conditions of \cref{lm:extremalD}, we can assign the leftover weight to triangles in $\mathcal{D}$ to obtain the desired fractional $\{C^{\stri}_{\ell},C^{\splus}_{\ell}\}$-decomposition. In practice, this means we need to verify that $\mathcal{D}$ satisfies the Hall-type condition of \cref{L:hallCondition}, that is, we need to show that for any set $X\subseteq V(G)$, we have $\frac{\ell-3}{3}\|\D_X\|-\|L[X]\| \geq 0$. We consider three cases depending on the size of $X$. 

First, we use that cross edges have significant weight in $\cD$ (\cref{lm:extremalD}\cref{D:crossmin}) to show that \cref{E:hallCondition} holds whenever $|X|$ is not too large. The intuition behind this is the following. Let $X\subseteq V(G)$ and suppose $|\overline{X}|$ is not too small. We want to show that $\|L[X]\|$ is relatively large compared to $\|\cD_X\|$. Since side edges are almost entirely covered by $\cD$ (\cref{lm:extremalD}\cref{D:side1,D:sideU}) and almost all of the cross edges accrue a non-trivial amount of weight in $\D$ (\cref{lm:extremalD}\cref{D:crossmin}), we have that the leftover $\|L[X]\|$ within $X$ cannot be excessively large compared to $\|\cD[X]\|$, the total weights of triangles which are fully within $X$.
But, as $|\overline{X}|$ is not too small by assumption, there are a significant number of cross edges between $X$ and $\overline{X}$ and most of these edges receive a non-trivial amount of weight in $\D$ (\cref{lm:extremalD}\cref{D:crossmin}). Thus, we obtain that $\|\cD_X\setminus \cD[X]\|$, the total weights of triangles incident to, but not fully within $X$, is non-trivial and can compensate any excess leftover within $\|L[X]\|$.

\begin{lemma}\label{lm:extremalsmall}
    Let $\ell \geq 73$, let $G$ be a near $\ell$-extremal graph, and suppose that $\cD$ is a fractional $K_3$-packing of $G$ satisfying \cref{lm:extremalD}\cref{D:cross,D:crossmax,D:crossmin,D:side1,D:sideU,D:Tsmall}. If $X\subseteq V(G)$ with $|X| \leq \frac{2\ell-9}{2\ell-5}n$, then $\frac{\ell-3}{3}\|\D_X\|-\|L[X]\| \geq 0$.
\end{lemma}

\begin{proof}
    We prove the contrapositive. Namely, suppose that $X\subseteq V(G)$ is such that $\|L[X]\|-\frac{\ell-3}{3}\|\D_X\|>0$. It suffices to show that $|\overline{X}| < \frac{4}{2\ell-5}n$.

    Observe that at most $|U||U \cap X|$ edges of $G[U]$ can be incident to a vertex in $X$ and that at most $\frac{1}{2}|U \cap X|^2$ of these can be in $G[U \cap X]$. By \cref{lm:extremalD}\cref{D:cross}, every triangle in $\cD$ contains exactly two edges of $\Gc$. Together with \cref{lm:extremalD}\cref{D:crossmin}, we obtain on the
    one hand that
    \begin{equation}\label{E:extremalsmallD}
        \|\cD_X\|\geq \mfrac{3}{2}\cdot \mfrac{2}{\ell}\bigl(e_G(X\cap A, X\cap B)+e_G(X\cap A, B\setminus X)+e_G(A\setminus X, X\cap B)-|U||U \cap X|\bigr).
    \end{equation}
    On the other hand, \cref{lm:extremalD}\cref{D:crossmin,D:side1,D:sideU} imply that
    \begin{align}
        \|L[X]\|&\leq \sum_{e \in E(\Gs[X])\setminus E(G_s[U])}w_L(e)+\mfrac{\ell-2}{\ell}\left(e_G(X\cap A, X\cap B)-\mfrac{1}{2}|U \cap X|^2\right) +\mfrac{1}{2}|U \cap X|^2 \notag\\
        &\leq \mfrac{1}{3}|U||X|+\mfrac{\ell-2}{\ell}e_G(X\cap A, X\cap B) +\mfrac{1}{\ell}|U \cap X|^2. \label{E:extremalsmallL}
    \end{align}
    We have $\|L[X]\|-\frac{\ell-3}{3}\|\D_X\| > 0$ by assumption and thus, combining \cref{E:extremalsmallD,E:extremalsmallL} and using $|U \cap X|^2 \leq |U||U \cap X|$, we obtain
    \begin{equation}\label{eq:Xsmall}
        0 < 
        \mfrac{1}{3}|U||X|+\mfrac{1}{\ell}e_G(X\cap A, X\cap B)+\mfrac{\ell-2}{\ell}|U||U \cap X|-\mfrac{\ell-3}{\ell}\bigl(e_G(X\cap A, B\setminus X)+e_G(A\setminus X, X\cap B)\bigr).
    \end{equation}
    We first obtain a crude lower bound on $|X|$. Using $\frac{\ell-2}{\ell}|U||U \cap X| \leq |U||X|$, we obtain
    \begin{align*}
        0&\stackrel{\text{\eqmakebox[crude]{\cref{eq:cU}}}}{<}  \mfrac{4}{3}|U||X|+\mfrac{1}{4\ell}|X|^2-\mfrac{\ell-3}{\ell}\left(|X|\cdot \mfrac{\delta(G)}{2}-2\cdot\mfrac{|X|^2}{4}\right)\\
        &\stackrel{\text{\eqmakebox[crude]{\cref{extremalDeg}}}}{\leq }|X|\left(\mfrac{2\ell-5}{4\ell}|X|-\mfrac{\ell-3}{4(\ell-1)}n+\mfrac{4}{3}|U|\right)\\
        &\stackrel{\text{\eqmakebox[crude]{}}}{\leq}|X|\left(\mfrac{1}{2}|X|-\mfrac{35}{144}n+\mfrac{4}{3}|U|\right)
    \end{align*}
    and so, using \cref{L:mAndU}\cref{uBound2},
    \begin{equation}\label{eq:Xsmallrough}
        |X|>\mfrac{35}{72}n-\mfrac{8}{3}|U|\geq \mfrac{n}{3}.
    \end{equation}
    With our crude bound in hand, we now work more carefully to prove the required upper bound on $|\overline{X}|$. From \cref{eq:Xsmall} we have
    \begin{align*}
        0 &< \mfrac{1}{3}|U||X|+\mfrac{1}{\ell}e_G(X\cap A, X\cap B)+\mfrac{\ell-2}{\ell}|U|^2-\mfrac{\ell-3}{\ell}\bigl(e_G(X\cap A, B\setminus X)+e_G(A\setminus X, X\cap B)\bigr)\\
        &\leq \mfrac{1}{3}|U||X|+\mfrac{1}{\ell}|X\cap A||X\cap B|+\mfrac{\ell-2}{\ell}|U|^2\\
        &\qquad -\mfrac{\ell-3}{\ell}\left(|X\cap A|\left(\mfrac{n}{2}-m-|X\cap B|\right)+\left(\mfrac{n}{2}+m-|X\cap A|\right)|X\cap B|-e_{\overline{G}}(A,B)\right)
        \\
        &= \mfrac{1}{3}|U||X|+\mfrac{1}{\ell}|X\cap A|\bigl(|X|-|X\cap A|\bigr)+\mfrac{\ell-2}{\ell}|U|^2\\
        &\qquad -\mfrac{\ell-3}{\ell}\left(\left(\mfrac{n}{2}+m\right)|X|-2|X\cap A|\bigl(|X|+m-|X\cap A|\bigr)-\Bigl(\mfrac{n^2}{4}-m^2-e(\Gc)\Bigr)\right)\\
        &= \mfrac{1}{3}|U||X|+\mfrac{\ell-2}{\ell}|U|^2+\mfrac{\ell-3}{\ell}\left(\mfrac{n^2}{4}-m^2-e(\Gc)\right)-\mfrac{\ell-3}{\ell}\left(\mfrac{n}{2}+m\right)|X|\\
        &\qquad -\mfrac{2\ell-5}{\ell}|X\cap A|^2+\mfrac{1}{\ell}\bigl((2\ell-5)|X|+2m(\ell-3)\bigr)|X\cap A|.
    \end{align*}
   The last expression is a concave quadratic in $|X\cap A|$, which takes its maximum at $|X\cap A|=\frac{|X|}{2}+\frac{\ell-3}{2\ell-5}m$.
   Thus,
   \begin{align*}
       0&\stackrel{\text{\eqmakebox[XandXbar]{}}}{<} \mfrac{1}{3}|U||X|+\mfrac{\ell-2}{\ell}|U|^2+\mfrac{\ell-3}{\ell}\left(\mfrac{n^2}{4}-m^2-e(\Gc)\right)+\mfrac{2\ell-5}{4\ell}|X|^2-\mfrac{\ell-3}{2\ell}n|X|+\mfrac{(\ell-3)^2}{\ell(2\ell-5)}m^2\\
       &\stackrel{\text{\eqmakebox[XandXbar]{}}}{\leq} \mfrac{1}{3}|U||X|+\mfrac{\ell-2}{\ell}|U|^2+\mfrac{\ell-3}{\ell}\left(\mfrac{n^2}{4}-e(\Gc)\right)+\mfrac{n}{4\ell}|X|-\mfrac{2\ell-5}{4\ell}|\overline{X}||X|-\mfrac{(\ell-3)(\ell-2)}{\ell(2\ell-5)}m^2\\
       &\stackrel{\text{\eqmakebox[XandXbar]{\cref{extremalCut}}}}{\leq} \mfrac{1}{3}|U||X|+\mfrac{\ell-2}{\ell}|U|^2+\mfrac{n^2}{\ell(\ell-1)}+\mfrac{n}{4\ell}|X|-\mfrac{2\ell-5}{4\ell}|\overline{X}||X|.
   \end{align*}
   Rearranging, we obtain
   \begin{align*}
       |\overline{X}|&\stackrel{\text{\eqmakebox[Xbar]{}}}{<}\mfrac{1}{2\ell-5}\left(n+\mfrac{4\ell}{3}|U|+\mfrac{1}{|X|}\left(\mfrac{4n^2}{\ell-1}+4(\ell-2)|U|^2\right)\right)\\
       &\stackrel{\text{\eqmakebox[Xbar]{\cref{eq:Xsmallrough}}}}{\leq} \mfrac{1}{2\ell-5}\left(n+\mfrac{4\ell}{3}|U|+\mfrac{12n}{\ell-1}+\mfrac{12(\ell-2)}{n}|U|^2\right)\\
       &\stackrel{\text{\eqmakebox[Xbar]{\cref{L:mAndU}\cref{uBound}}}}{\leq} \mfrac{n}{2\ell-5}\left(1+\mfrac{140\ell}{(\ell-3)(\ell-1)}+\mfrac{12}{\ell-1}+\mfrac{132300(\ell-2)}{(\ell-3)^2(\ell-1)^2}\right)\\
       &\stackrel{\text{\eqmakebox[Xbar]{$\ell \geq 73$}}}{<} \mfrac{4n}{2\ell-5},
   \end{align*}
   as desired.
\end{proof}

Next, we use that each individual triangle has relatively small weight in $\cD$ (\cref{lm:extremalD}\cref{D:Tsmall}) \and that triangles within $U$ have no weight at all in $\cD$ (\cref{lm:extremalD}\cref{D:cross}) to show that if the entire vertex set satisfies \cref{E:hallCondition}, then so does any very large $X\subseteq V(G)$. The intuition behind this is the following. 
Suppose that $X\subseteq V(G)$ is a largest set such that \cref{E:hallCondition} does not hold and suppose for a contradiction that $X\neq V(G)$.
If $\overline{X}\subseteq U$, then $\|\cD\|=\|\cD_X\|$ as triangles in $U$ receive no weight at all in $\cD$ (\cref{lm:extremalD}\cref{D:cross}),
and so \cref{E:hallCondition} must also fail for the entire vertex set $V(G)$ as $\|L\|\geq \|L[X]\|$, a contradiction. Therefore, there exists $z\in \overline{X}\setminus U$. We claim that $X\cup \{z\}$ also fails \cref{E:hallCondition}, which contradicts the maximality of $X$.
Indeed, $\|\cD_{X\cup\{z\}}\|$ is obtained from $\|\cD_X\|$ by adding the weights of the triangles of $\cD$ which are incident to $z$ and fully within $\overline{X}$. However, $|\overline{X}|$ is very small by assumption, so as no triangle is assigned a lot of weight (\cref{lm:extremalD}\cref{D:Tsmall}), we have on the one hand that $\|\D_{X\cup\{z\}}\|\approx \|\cD_X\|$. On the other hand, 
the difference between $\|L[X\cup\{z\}\|$ and $\|L[X]\|$ is at least the weight of $L$ on the edges between $X$ and $z$, and this is comparatively large since $X$ is very large by assumption and no cross edge in $\overline{U}$ can accumulate too much weight in $\D$ (\cref{lm:extremalD}\cref{D:crossmax}).
Thus, \cref{E:hallCondition} must also fail for $X\cup \{z\}$, a contradiction to the maximality of $X$.

\begin{lemma}\label{lm:extremalXbar}
    Let $\ell \geq 73$, let $G$ be a near $\ell$-extremal graph, and suppose that $\cD$ is a fractional $K_3$-packing of $G$ satisfying \cref{lm:extremalD}\cref{D:cross,D:crossmax,D:crossmin,D:side1,D:sideU,D:Tsmall}. If $X\subseteq V(G)$ such that $|X| > \frac{2\ell-9}{2\ell-5}n$, then $\frac{\ell-3}{3}\|\cD_X\|-\|L[X]\| \geq \frac{\ell-3}{3}\|\cD\|-\|L\|$.
\end{lemma}

\begin{proof}
Of the subsets of $V(G)$ of cardinality greater than $\frac{2\ell-9}{2\ell-5}n$, let $X$ be one such that $\frac{\ell-3}{3}\|\cD_X\|-\|L[X]\|$ is minimal and, subject to this, $|X|$ is maximal. It suffices to show that $X=V(G)$.

If $\overline{X} \subseteq U$, then $w_{\D}(T)=0$ for each $T \in \T(G[\overline{X}])$ by \cref{lm:extremalD}\cref{D:cross}, which implies that
\[\mfrac{\ell-3}{3}\|\cD_X\|-\|L[X]\|=\mfrac{\ell-3}{3}\|\cD\|-\|L[X]\|\geq \mfrac{\ell-3}{3}\|\cD\|-\|L\|,\]
and so $X=V(G)$ by definition, as desired.

We may therefore assume that there exists a vertex $z \in \overline{X} \setminus U$. We will find a contradiction. Let $\D'_z$ be the fractional $K_3$-packing induced from $\D$ by the triangles that contain $z$ and two other vertices in $\overline{X}$.
By definition of $X$, we have
\begin{align}
    0 &< \left(\mfrac{\ell-3}{3}\|\cD_{X\cup\{z\}}\|-\|L[X\cup\{z\}]\|\right)-\left(\mfrac{\ell-3}{3}\|\cD_X\|-\|L[X]\|\right)\nonumber\\
    &=\mfrac{\ell-3}{3}\|\D'_z\|-\sum_{y \in N_G(z) \cap X}\bigl(1-w_\cD(yz)\bigr).\label{eq:largeX}
\end{align}
But using \cref{lm:extremalD}\cref{D:crossmax}, \cref{eq:cgood}, and \cref{L:mAndU}\cref{uBound2}, we have on the one hand that
\begin{equation}\label{E:LzBound}
    \sum_{y \in N_G(z) \cap X}\bigl(1-w_\cD(yz)\bigr) \geq \mfrac{3}{5}\bigl|(N_{\Gc}(z) \cap X) \setminus U\bigr| \geq \mfrac{3}{5}\left(\mfrac{9n}{20}-|\overline{X}|-|U|\right) > \mfrac{n}{5}.
\end{equation}
On the other hand, there is at most one triangle in $\D'_z$ for each edge in $G[\overline{X}\setminus \{z\}]$ and hence there are at most $\frac{1}{2}|\overline{X}|^2$ such triangles. Moreover, \cref{lm:extremalD}\cref{D:Tsmall} implies that $\|T\| \leq 3\cdot\frac{4}{n}=\frac{12}{n}$ for each $T \in \D'_z$, and so
\begin{equation}\label{E:DzBound}
    \mfrac{\ell-3}{3}\|\D'_z\| \leq \mfrac{\ell-3}{3} \cdot \mfrac{12}{n} \cdot \mfrac{1}{2}|\overline{X}|^2 \leq \mfrac{2(\ell-3)}{n} \cdot \left(\mfrac{4n}{2\ell-5}\right)^2 < \mfrac{16n}{2\ell-5} < \mfrac{n}{5}.
\end{equation}
Together, \cref{E:LzBound} and \cref{E:DzBound} contradict \cref{eq:largeX}.
This concludes the proof.
\end{proof}

Finally, we use that side edges are extremely well covered by $\mathcal{D}$ (\cref{lm:extremalD}\cref{D:side1,D:sideU,D:cross}) to show that the entire vertex set satisfies \cref{E:hallCondition}, that is, that $\|\D\|\geq \frac{3}{\ell}e(G)$.
The intuition behind this is the following. If $G\in \G_\ell$, in which case 
$|A|=|B|$ and there are no exceptional vertices, then precisely $\frac{1}{\ell}$ of the edges of $G$ are side edges and so, as these are fully covered by $\cD$ (\cref{lm:extremalD}\cref{D:side1}) and we ignore side triangles (\cref{lm:extremalD}\cref{D:cross}), we automatically obtain that $\|\D\|$ is sufficiently large. When $G$ has more edges than the graphs in $G_\ell$ and/or $|A|\neq|B|$ (i.e.\ $m>0$), the proportion of side edges is even larger, so $\|\D\|$ is even larger than desired. The only issue is the exceptional vertices in $U$, as side edges incident to $U$ may not be fully covered by $\cD$. However, recall that exceptional vertices have relatively high side degree by definition, so the presence of such vertices increases the proportion of side edges in $G$ significantly. Moreover, we know that there are relatively few vertices in $U$ (\cref{L:mAndU}), so most of those extra side edges are between exceptional and non-exceptional vertices. As these edges still have large weight (\cref{lm:extremalD}\cref{D:sideU}), the extra proportion of side edges counterbalance the small loss of weight on exceptional side edges. That is, perhaps counter-intuitively, exceptional vertices make it easier to satisfy \cref{E:hallCondition}.

\begin{lemma}\label{lm:extremalfull}
    Let $\ell \geq 73$, let $G$ be a near $\ell$-extremal graph, and suppose that $\cD$ is a fractional $K_3$-packing of $G$ satisfying \cref{lm:extremalD}\cref{D:cross,D:crossmax,D:crossmin,D:side1,D:sideU,D:Tsmall}. Then $\|\cD\| \geq \frac{3}{\ell}e(G)$ and hence $\frac{\ell-3}{3}\|\cD\|-\|L\| \geq 0$.
\end{lemma}

\begin{proof}
If $\|\cD\| \geq \frac{3}{\ell}e(G)$, then $\frac{\ell-3}{3}\|\cD\|-\|L\| \geq 0$ follows immediately since $\|L\|=e(G)-\|\D\|$. By \cref{lm:extremalD}\cref{D:cross}, we have $\|\cD\|=3\sum_{e\in E(\Gs)}w_{\D}(e)$. Therefore, to prove that $\|\cD\| \geq \frac{3}{\ell}e(G)$, it suffices to show that
    \begin{align}\label{eq:full}
        \mfrac{2}{3}\left(\|\D\|-\mfrac{3}{\ell}e(G)\right)=2\sum_{e\in E(\Gs)}w_{\D}(e)-\mfrac{1}{\ell}\sum_{v\in V(G)}d_G(v)
    \end{align}
    is nonnegative.
    
    Let $E_U$ denote the set of all edges of $G$ that are incident to at least one vertex in $U$. By \cref{lm:extremalD}\cref{D:side1,D:sideU}, we have
    \begin{align}
        2\sum_{e\in E(\Gs)}w_{\D}(e)&=2e(\Gs)-2\sum_{uv \in E_U}\bigl(1-w_{\D}(uv)\bigr) \notag\\
        &= \sum_{v\in V(G)} d_{\Gs}(v)-2\sum_{uv \in E_U}\bigl(1-w_{\D}(uv)\bigr) \notag\\
        &\geq \sum_{v\in V(G)} d_{\Gs}(v)-2\left(\sum_{u\in U}d_{\Gs}(u)-|U|^2\right)\left(1-\mfrac{2}{3}\right)-|U|^2 \notag\\
        &= \sum_{v\in V(G)} d_{\Gs}(v)-\mfrac{2}{3}\sum_{u\in U}d_{\Gs}(u)-\mfrac{1}{3}|U|^2 \notag\\
        &= \sum_{v\in V(G) \setminus U} d_{\Gs}(v)+\mfrac{1}{3}\sum_{u\in U}d_{\Gs}(u)-\mfrac{1}{3}|U|^2. \label{E:fullsideedges}
    \end{align}
    Since $d_{\Gs}(v)=d_{G}(v)-d_{\Gc}(v)$ for all $v \in V(G)$, we have
    \begin{align}
        \sum_{v\in V(G)\setminus U} d_{\Gs}(v)-\mfrac{1}{\ell}\sum_{v\in V(G)\setminus U}d_G(v)
        &\stackrel{\text{\eqmakebox[fullsumdiff]{}}}{=} \mfrac{\ell-1}{\ell} \sum_{v\in V(G)\setminus U} d_{G}(v)-\sum_{v\in V(G)\setminus U}d_{\Gc}(v) \notag\\
        &\stackrel{\text{\eqmakebox[fullsumdiff]{\cref{extremalDeg}}}}{\geq} \mfrac{\ell-1}{\ell}\cdot \mfrac{\ell n}{2(\ell-1)}(n-|U|)-\left(2|A||B|-|U|\left(\mfrac{n}{2}-m\right)\right) \notag\\
        &\stackrel{\text{\eqmakebox[fullsumdiff]{}}}{=}2m^2-m|U|. \label{E:fullsumdiff}
    \end{align}
    Combining these facts, we obtain 
    \begin{align}
        \mfrac{2}{3}\left(\|\D\|-\mfrac{3}{\ell}e(G)\right)
        &\stackrel{\text{\eqmakebox[extremalfullmquadratic]{\cref{eq:full},\cref{E:fullsideedges}}}}{\geq} \sum_{v\in V(G) \setminus U} d_{\Gs}(v)+\mfrac{1}{3}\sum_{u\in U}d_{\Gs}(u)-\mfrac{1}{3}|U|^2 -\mfrac{1}{\ell}\sum_{v\in V(G)}d_G(v) \notag\\
        &\stackrel{\text{\eqmakebox[extremalfullmquadratic]{\cref{E:fullsumdiff}}}}{\geq} 2m^2-m|U|+\left(\mfrac{1}{3}-\mfrac{1}{\ell}\right)\sum_{u\in U}d_{\Gs}(u)-\mfrac{1}{3}|U|^2-\mfrac{1}{\ell}\sum_{u\in U}d_{\Gc}(u)\notag\\
        &\stackrel{\text{\eqmakebox[extremalfullmquadratic]{}}}{\geq} 2m^2-m|U|+\left(\mfrac{1}{3}-\mfrac{1}{73}\right)\mfrac{n}{20}|U|-\mfrac{1}{3}|U|^2-\mfrac{1}{73}|U|\left(\mfrac{n}{2}+m\right)\notag\\
        &\stackrel{\text{\eqmakebox[extremalfullmquadratic]{}}}{=} 2m^2-\mfrac{74}{73}m|U|+\mfrac{2n}{219}|U|-\mfrac{1}{3}|U|^2. \label{E:extremalfullmquadratic}
    \end{align}
    This is a convex quadratic in $m$ that takes its minimum at $m=\frac{37}{146}|U|$. Thus,
    \begin{equation}\label{E:extremalfullfinal}
        \mfrac{2}{3}\left(\|\D\|-\mfrac{3}{\ell}e(G)\right)
        \geq \mfrac{1}{31974}|U|\bigl(292n - 14765|U|\bigr) \stackrel{\text{\cref{L:mAndU}\cref{uBound2}}}{\geq}0,
    \end{equation}
    as desired.
\end{proof}

\begin{proof}[Proof of \cref{lm:extremal}]
Let $\D$ be a fractional $K_3$-packing of $G$ satisfying \cref{lm:extremalD}\cref{D:cross,D:crossmax,D:crossmin,D:side1,D:sideU,D:Tsmall} and let $L$ be the leftover of $\D$.
It follows from \cref{lm:extremalsmall,lm:extremalXbar,lm:extremalfull} that $\frac{\ell-3}{3}\|\D_X\|-\|L[X]\| \geq 0$
holds for all $X \subseteq V(G)$. Thus the result follows by \cref{L:hallCondition}.
\end{proof}

\section{Proof of Theorem \ref{T:main} and concluding remarks}\label{sec:conclusion}

\subsection{Odd cycles of lengths at least 73}

We are now ready to determine the decomposition threshold for cycles $C_\ell$ of odd length $\ell\geq 73$.

\begin{proof}[Proof of \cref{T:main}]
As we saw in the introduction, $\delta_{C_\ell} \geq \frac{1}{2}+\frac{1}{2(\ell-1)}$.  From \cref{lm:nonextremal,lm:extremal} we have that every $n$-vertex graph $G$ with minimum degree at least $\frac{1}{2}+\frac{1}{2(\ell-1)}$ has a fractional $\{C^{\stri}_{\ell},C^{\splus}_{\ell}\}$-decomposition. Thus, since $C^{\stri}_{\ell}$ and $C^{\splus}_{\ell}$ are condensations of $C_\ell$ by \cref{lm:cond}, we have that $\delta^{0+}_{C_\ell} \leq \frac{1}{2}+\frac{1}{2(\ell-1)}$ by \cref{T:condensations}. So $\delta_{C_\ell} \leq \frac{1}{2}+\frac{1}{2(\ell-1)}$ by \cref{T:reduction}.
\end{proof}

\subsection{Smaller odd cycles}
The value of $73$ is clearly an artefact of our proof and we believe that \cref{conj:odd} also holds for all odd $5\leq \ell \leq 71$.
Also note that the value of $73$ is not optimised. Indeed, some parameters in the extremal case (e.g.\ the threshold of $\frac{n}{20}$ for a vertex to be exceptional) are not fully optimised, and we simply pick sensible values that keep all the calculations relatively neat. There are also various places in the extremal case where one could refine our arguments. However, we opt against doing so as this will not lead to a proof for $\ell=5$ (say). Indeed, the main barrier of our method is the threshold between our extremal and non-extremal cases, which measures how far a graph is from those in $\G_\ell$ as a function of $\frac{1}{(\ell-3)(\ell-1)}$. When $\ell$ is not too small, $\frac{1}{(\ell-3)(\ell-1)}$ is quite small, and so near $\ell$-extremal graphs are indeed very close to the extremal graphs in $\G_\ell$. But, for instance, near $5$-extremal graphs may have up to $\frac{n^2}{16}$ extra edges in total and up to $\frac{n^2}{8}$ fewer edges in their maximum edge cut, and so may not actually look like a graph in $\G_5$.

To cover smaller values of $\ell$, the crux therefore lies in improving the non-extremal case. One option could be to improve our analysis to show that the simple fractional $K_3$-packing given by \cref{L:decentPartialDecomp} actually works for more graphs. This would improve the threshold at which a graph becomes near $\ell$-extremal and our arguments for the extremal case should then work mutatis mutandis for smaller values of $\ell$. 

However, it is very likely that, for small values of $\ell$, the simple packing $\cD$ from \cref{L:decentPartialDecomp} is not sufficient. This means one would need a more refined way of constructing the packing when $G$ is not too close to the extremal graphs in $\G_\ell$ but some structure is still discernable. There are natural candidates for such a suitable packing, however the difficulty lies in finding a good analysis of its properties to be able to prove the Hall-type condition \cref{E:hallCondition}.

\subsection{Other graphs}
Observe that what we showed is in fact that any graph of minimum degree at least $(\frac{1}{2}+\frac{1}{2(\ell-1)})n$ has a fractional $\{C^{\stri}_{\ell},C^{\splus}_{\ell}\}$-decomposition. Together with \cref{T:condensations,T:reduction}, this gives rise to the desired $C_\ell$-decomposition.
But note that there are many non-bipartite graphs other than $C_\ell$ which have both $C^{\stri}_{\ell}$ and $C^{\splus}_{\ell}$ as condensations. For example, this is true of any unicyclic graph with $\ell$ edges in which the cycle has odd length and contains two adjacent vertices of degree 2.  Together with \cref{T:condensations}, our proof therefore implies that $\delta_F^{0+}=\frac{1}{2}+\frac{1}{2(\ell-1)}$ for any such graph $F$. However, we cannot determine $\delta_F$ for such graphs because we lack an analogue to \cref{T:reduction}.

\section*{Acknowledgments}

We thank Peter Dukes, Barbara Maenhaut, and Darryn Bryant for helpful discussions.

\bibliography{refs}
\bibliographystyle{abbrv}

\end{document}